\documentclass{article}

\usepackage[margin=2cm]{geometry}
\usepackage{amsmath,amssymb}
\usepackage{amsmath}
\usepackage{amsthm}
\usepackage{enumitem}
\usepackage{multicol}
\usepackage{graphicx}
\usepackage{url}
\usepackage[colorlinks]{hyperref}

\usepackage{color}

\usepackage{listings}
\usepackage{tikz}
\tikzset{
every node/.style={circle, draw, inner sep=2pt},
every picture/.style={thick}
}
\usetikzlibrary{matrix,arrows,calc}

\newtheorem{theorem}{Theorem}
\newtheorem{lemma}[theorem]{Lemma}
\newtheorem{proposition}[theorem]{Proposition}
\newtheorem{corollary}[theorem]{Corollary}

\theoremstyle{definition}
\newtheorem{definition}[theorem]{Definition}
\newtheorem{observation}[theorem]{Observation}
\newtheorem{remark}[theorem]{Remark}
\newtheorem{example}[theorem]{Example}
\newtheorem{problem}{Problem}
\newtheorem{question}[theorem]{Question}
\newtheorem{conjecture}[theorem]{Conjecture}

\newcommand{\diag}{\operatorname{diag}}
\newcommand{\minors}{\operatorname{minors}}

\DeclareMathOperator{\trs}{trs}
\DeclareMathOperator{\dist}{dist}
\DeclareMathOperator{\diam}{diam}
\DeclareMathOperator{\snf}{SNF}

\usepackage[textwidth=2.5cm, textsize=small, colorinlistoftodos]
{todonotes}

\newcommand{\comment}[1]{}

\begin{document}


\title{Distance ideals of digraphs}


\author{Carlos A. Alfaro$^a$, Teresa I. Hoekstra-Mendoza$^b$, Juan Pablo Serrano$^c$, Ralihe R. Villagrán$^d$
\\ \\
{\small $^a$ Banco de M\'exico} \\
{\small Mexico City, Mexico}\\
{\small {\tt carlos.alfaro@banxico.org.mx}}
\\
{\small $^b$ Centro de Investigaci\'on en Matem\'aticas}\\
{\small Guajanajuato, Mexico}\\
{\small\tt maria.idskjen@cimat.mx}
\\
{\small $^c$Departamento de Matem\'aticas}\\
{\small Centro de Investigaci\'on y de Estudios Avanzados del IPN}\\
{\small Apartado Postal 14--740 }\\
{\small 07000 Mexico City, Mexico. } \\
{\small {\tt jpserranop@math.cinvestav.mx}}
\\
{\small $^d$Department of Mathematical Sciences} \\
{\small Worcester Polytechnic Institute, Worcester, USA}\\
{\small {\tt
rvillagran@wpi.edu}} \\
\\
}
\date{}


\maketitle

\begin{abstract}
    We focus on strongly connected, strong for short, digraphs since in this setting distance is defined for every pair of vertices. 
    Distance ideals generalize the spectrum and Smith normal form of several distance matrices associated with strong digraphs.
    We introduce the concept of pattern which allow us to characterize the family $\Gamma_1$ of digraphs with only one trivial distance ideal over ${\mathbb Z}$.
    This result generalizes an analogous result for undirected graphs that states that connected graphs with one trivial ideal over $\mathbb{Z}$ consists of either complete graphs or complete bipartite graphs.
    It turns out that the strong digraphs in $\Gamma_1$ consists in the circuit with 3 vertices and a family $\Lambda$ of strong digraphs that contains complete graphs and complete bipartite graphs, regarded as digraphs.
    We also compute all distance ideals of some strong digraphs in the family $\Lambda$.
    Then, we explore circuits, which turn out to be an infinite family of minimal forbidden digraphs, as induced subdigraphs, for the strong digraphs in $\Gamma_1$.
    This suggests that a characterization of $\Gamma_1$ in terms of forbidden induced subdigraphs is harder than using patterns.  

    
\end{abstract}


\section{Introduction}

A \emph{digraph} $G$ is a pair $(V,A)$ of vertices and arcs.
In this study, we focus on simple digraph, that is, multiple arcs and loops are not allowed.
A \emph{walk} $W$ in a digraph $G$ is an alternating sequence $u_1a_1u_2 \cdots u_{k-1}a_{k-1}u_k$ of vertices and arcs such that the tail of arc $a_i=(u_iu_{i+1})$ is $u_i$ and the head is $u_{i+1}$ for every $i\in[k-1]$, where $[k]:=\{1,..., k\}$.
A walk starting in vertex $u$ and ending in vertex $v$ is called $uv$-\emph{walk}.
A $uv$-\emph{path} of length $k-1$ is a sequence $u_1u_2\cdots u_k$ of distinct (except for the start and the end) vertices such that $u=u_1$, $v=u_k$, and for each $i\in\{1,\dots,k-1\}$, there exists an arrow $(u_iu_{i+1})\in A(G)$. 
The \emph{circuit} $\overrightarrow{C_n}$ is a closed path, that is $u=v$, of length $n$.
A digraph $G$ is \emph{strongly connected} or just \emph{strong} if for every pair of distinct vertices $u$ and $v$ in $G$, then there exists a $uv$-walk and a $vu$-walk\footnote{At \url{https://oeis.org/A035512} can be found the first numbers in the sequence of the number strong digraphs with $n$ vertices.}.
Circuits are basic structures of strong digraphs, since circuits naturally appear as subdigraphs.
Note that any strong digraph with at least 2 vertices without multiple arcs contains at least one circuit as an induced subdigraph.
This is not true if multiple arcs are allowed.
In the following, digraphs are considered strong unless otherwise stated, the reason is that we will need that for any pair $u,v$ of vertices there exists at least one $uv$-walk, then we will be able to define the distance from $u$ to $v$.
The {\it distance} $dist_G(u,v)$ from vertex $u$ to vertex $v$ is the number of arcs in a shortest $uv$-walk. 
The {\it distance matrix} $D(G)$ of $G$ is the matrix with rows and columns indexed by the vertices of $G$ with the $uv$-entry equal to $\dist_G(u,v)$.
Distance matrices were introduced in \cite{GP} by Graham and Pollack regarding the study of a data communication problem.
This problem involved finding appropriate addresses for a message to move efficiently through a series of loops from its origin to the destination, choosing the best route at each switching point.
The \emph{transmission} $\trs(u)$ of vertex $u$ is $\sum_{v\in V(G)}\dist(u,v)$, thus
$\trs(G)$ is the diagonal matrix with the transmissions of the vertices of $G$ in the diagonal, and let $\deg(G)$ denote the diagonal matrix with the out-degrees of the vertices of $G$ in the diagonal.
Then, we are able to define other closed related matrices: the \emph{distance Laplacian matrix} $D^L(G):=\trs(G)-D(G)$, the \emph{signless distance Laplacian matrix} $D^Q(G):=\trs(G)+D(G)$, the \emph{degree-distance matrix} $D^{\deg}(G):=\deg(G)-D(G)$, and the \emph{signless degree-distance matrix} $D^{\deg}_+(G):=\deg(G)+D(G)$.
Some properties of these matrices have been studied in \cite{degreedistance,ah2013}.

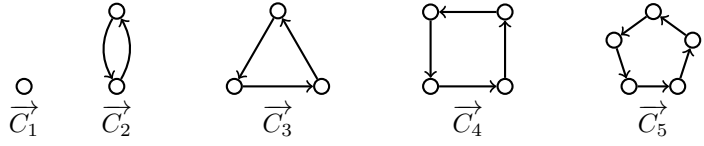
\begin{figure}
    \centering
    \begin{tabular}{c@{\extracolsep{0.8cm}}cccc}
        \begin{tikzpicture}[scale=2,thick]
		\tikzstyle{every node}=[minimum width=0pt, inner sep=2pt, circle]
			\draw (0,0) node[draw] (0) {};
		\end{tikzpicture}
        &
        \begin{tikzpicture}[scale=1,thick]
		\tikzstyle{every node}=[minimum width=0pt, inner sep=2pt, circle]
			\draw (0,0) node[draw] (0) {};
			\draw (0,1) node[draw] (1) {};
			\draw (0) edge[bend right,->] (1);
            \draw (1) edge[bend right,->] (0);
		\end{tikzpicture}
        &
        \begin{tikzpicture}[scale=0.67,thick]
		\tikzstyle{every node}=[minimum width=0pt, inner sep=2pt, circle]
			\draw (-30:1) node[draw] (0) {};
			\draw (90:1) node[draw] (1) {};
			\draw (210:1) node[draw] (2) {};
			\draw (0) edge[->] (1);
            \draw (1) edge[->] (2);
            \draw (2) edge[->] (0);
		\end{tikzpicture}
        &
        \begin{tikzpicture}[scale=0.7,thick]
		\tikzstyle{every node}=[minimum width=0pt, inner sep=2pt, circle]
			\draw (45:1) node[draw] (0) {};
			\draw (135:1) node[draw] (1) {};
			\draw (225:1) node[draw] (2) {};
			\draw (-45:1) node[draw] (3) {};
			\draw (0) edge[->] (1);
            \draw (1) edge[->] (2);
            \draw (2) edge[->] (3);
            \draw (3) edge[->] (0);
		\end{tikzpicture}
        &
        \begin{tikzpicture}[scale=0.55,thick]
		\tikzstyle{every node}=[minimum width=0pt, inner sep=2pt, circle]
			\foreach \x in {1, 2, 3, 4, 5}
    		{
    			\draw (18 + 72*\x:1) node[draw] (\x) {};
    		}
        	\draw (5) edge[->] (1);
            \draw (1) edge[->] (2);
            \draw (2) edge[->] (3);
            \draw (3) edge[->] (4);
            \draw (4) edge[->] (5);
		\end{tikzpicture}
		\\
		$\overrightarrow{C_1}$ & $\overrightarrow{C_2}$ & $\overrightarrow{C_3}$ & $\overrightarrow{C_4}$ & $\overrightarrow{C_5}$
    \end{tabular}
    \caption{The first five small circuits.}
    \label{fig:circuits}
\end{figure}

Given a strong digraph $G$ with $n$ vertices and a set of indeterminates $X_G=\{x_u \, : \, u\in V(G)\}$.
Whenever possible, we will use $X$ instead of $X_G$ to simplify the notation.
We define the matrix $D_X(G)$ as  $\diag(x_1, \dots,x_n)+D(G)$.
Note we can recover the distance matrix from $D_X$ by evaluating $X$ at the zero vector, that is, $D(G)=D_X(G)|_{X=\bf{0}}$.
Let $\mathfrak{R}[X]$ be the polynomial ring over a commutative ring $\mathfrak{R}$ in the variables $X$.
For all $i\in[n]$, the $i${\it-th distance ideal} $I^\mathfrak{R}_i(G)$ of $G$ is the ideal, over $\mathfrak{R}[X]$, given by
$\langle \minors_i(D_X(G))\rangle$,
where $n$ is the number of vertices of $G$ and ${\rm minors}_i(D_X(G))$ is the set of the determinants of the $i\times i$ submatrices of $D_X(G)$.
Distance ideals were defined in \cite{at} as a generalization of the Smith normal form of distance matrix and the distance spectra of graphs.
Often we will be interested in the case when $\mathfrak{R}$ is the ring $\mathbb{Z}$ of integer numbers.
In such cases we will omit $\mathfrak{R}$ in the notation of distance ideal, that is, we will use $I_i(G)$ to denote the $i$-{\it th} distance ideal of $G$ over $\mathbb{Z}[X]$.

\begin{example}
    Consider the circuit $\overrightarrow{C_3}$ with vertices $v_0$, $v_1$ and $v_2$.
    Thus
    \[
    D_X\left(\overrightarrow{C_3}\right)=
    \begin{bmatrix}
        x_0 &   1 &   2\\
          2 & x_1 &   1\\
          1 &   2 & x_2\\
    \end{bmatrix}.
    \]
    Then $I_1\left(\overrightarrow{C_3}\right)=\langle1\rangle$, and $I_2\left(\overrightarrow{C_3}\right)$ is generated by the Gröbner basis $\langle x_0 + 3, x_1 + 3, x_2 + 3, 7\rangle$, meanwhile $I_3\left(\overrightarrow{C_3}\right)$ is generated by the determinant of $D_X\left(\overrightarrow{C_3}\right)$, that is, $x_0x_1x_2 - 2x_0 - 2x_1 - 2x_2 + 9$.
    On the other hand, consider the complete digraph $K_3$.
    In this case
    \[
    D_X(K_3)=
    \begin{bmatrix}
        x_0 &   1 &   1\\
          1 & x_1 &   1\\
          1 &   1 & x_2\\
    \end{bmatrix}.
    \]
    Therefore $I_1(K_3)=\langle1\rangle$, and $I_2(K_3)$ is generated by the Gröbner basis $\langle x_0 - 1, x_1 - 1, x_2 - 1\rangle$, meanwhile $I_3(K_3)$ is generated by the determinant of $D_X(K_3)$, that is, $x_0x_1x_2 - x_0 - x_1 - x_2 + 2$.
\end{example}

Some special distance ideals are the {\it univariate distance ideals} which are defined from the matrix $D_t(G)=t{\sf I}_n+D(G)$, where ${\sf I}_n$ is the identity matrix of order $n$.
Then, the $i${\it-th univariate distance ideal} $U_i^{\mathcal{R}}(G)$ is the ideal $\langle\minors_i(D_t(G))\rangle$.
As previously mentioned, we will omit the base ring $\mathfrak R$ when it is $\mathbb Z$.
Thus $U_1(\overrightarrow{C_3})=\langle1\rangle$, $U_2(\overrightarrow{C_3})=\langle t + 3, 7\rangle$, and $U_3(\overrightarrow{C_3})=\langle t^3 - 6t + 9\rangle$.
On the other hand $U_1(K_3)=\langle1\rangle$, $U_2(K_3)=\langle t - 1\rangle$, and $U_3(K_3)=\langle t^3 - 3t + 2\rangle$.
Note, univariate distance ideals can be obtained by evaluating the variables at $x_i=t$, however they can also be computed directly.
In Appendix~\ref{appendix:distanceideals} we give some {\tt sagemath} codes that compute distance ideals and univariate distance ideals.

\begin{figure}[ht]
    \centering
    \begin{tabular}{c@{\extracolsep{2cm}}c}
         \includegraphics[width=4.5cm]{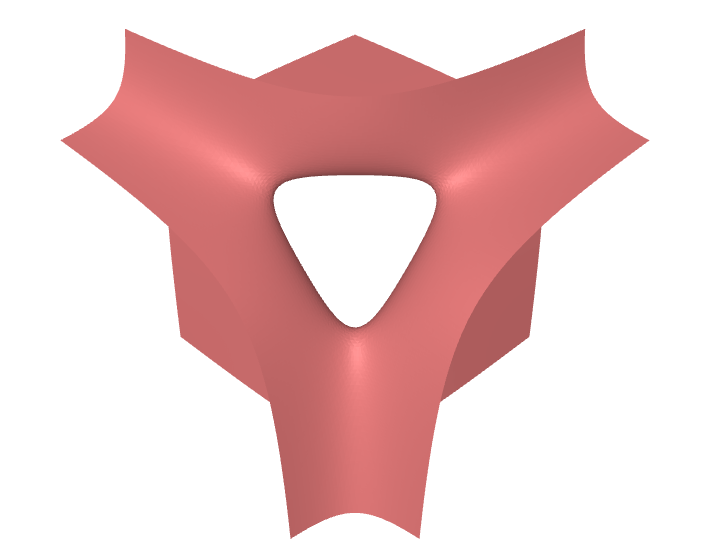}
         &
         \includegraphics[width=4cm]{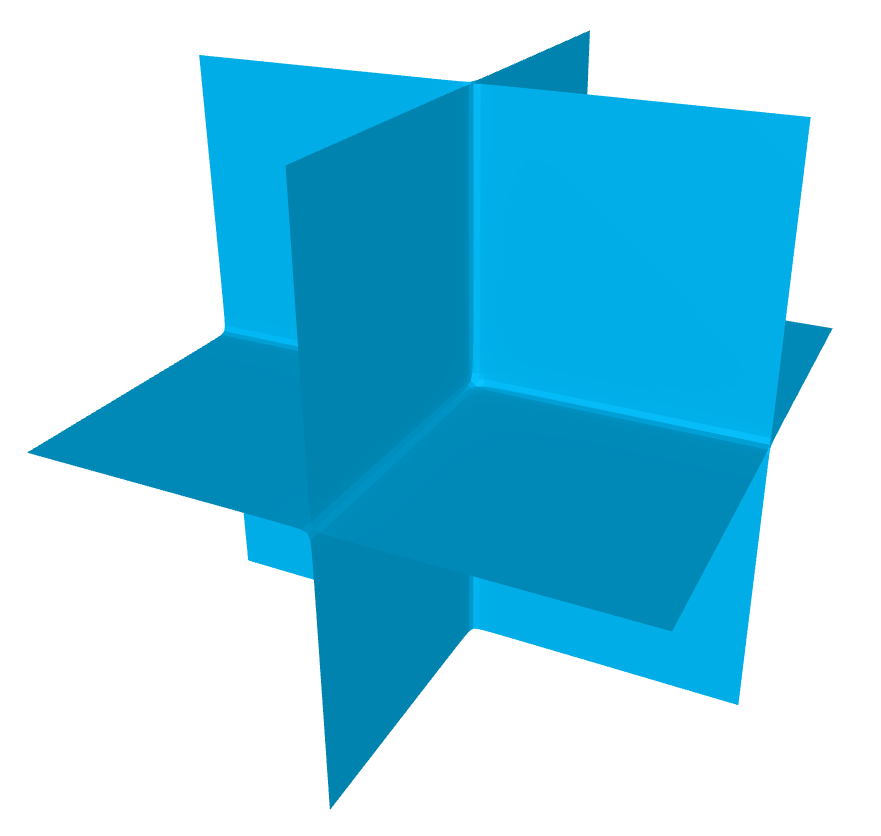}
    \end{tabular}
    \caption{Partial views of varieties in $\mathbb{R}^3$ associated with the third distance ideals of the circuit with 3 vertices (left) and the complete digraph with 3 vertices (right).}
    \label{fig:varieties}
\end{figure}

Let $I\subseteq \mathfrak{R}[X]$ be an ideal in $\mathfrak{R}[X]$, where $\mathfrak{R}$ is a commutative ring.
The \emph{variety} $V(I)$ of $I$ is defined as the set of common roots between polynomials in $I$.
In many contexts, it will be more convenient to consider an extension $\mathfrak{F}$ of $\mathfrak{R}$ to define the variety $V^\mathfrak{F}(I):=\left\{ {\bf a}\in \mathfrak{F}^l : f({\bf a}) = 0 \text{ for all } f\in I \right\}$.
Figure \ref{fig:varieties} shows partial views of varieties associated with the third distance ideals of some digraphs, and in Appendix~\ref{appendix:varieties}, it is the code to plot them.

The following result shows the contention with varieties, it was already observed in \cite{Nbook}, but it was specifically explored for distance ideals in \cite{akm}.
For any digraph $G$, it holds that
\[
\langle 1\rangle \supseteq I_1(G) \supseteq \cdots \supseteq I_n(G) \supseteq \langle 0\rangle
\]
and
\[
\emptyset\subseteq V^\mathfrak{F}(I_1(G)) \subseteq \cdots \subseteq V^\mathfrak{F}(I_n(G)) \subseteq \mathfrak{F}^l.
\]

Smith normal forms have been useful in understanding algebraic properties of combinatorial objects, see \cite{stanley}.
For instance, computing the Smith normal form of the adjacency or Laplacian matrix is a standard technique used to determine the Smith group and the critical group of a graph, see \cite{av,merino,rushanan}.
Little is known about the Smith normal forms of distance matrices.
In \cite{HW}, the Smith normal forms of the distance matrices were determined for trees, wheels, cycles, and complements of cycles and are partially determined for complete multipartite graphs.
In \cite{BK}, the Smith normal form was obtained for the distance matrices of
unicyclic graphs and wheel graphs with trees attached to each
vertex.
It is known that the Smith normal form may not exist in many commutative rings, for example consider the following matrix with entries in the ring $\mathbb{Z}[x]$
\[
\begin{bmatrix}
2 & 0\\
0 & x
\end{bmatrix}.
\]

Smith normal forms can be computed using row and column operations.
Two $n\times n$ matrices $M$ and $N$ are \emph{equivalent}, denoted by $M\sim N$, if there exist $P,Q\in GL_n(\mathbb{Z})$ such that $N=PMQ$.
That is, $M$ can be transformed to $N$ by applying the following elementary row and column operations which are invertible over the ring of integers:
\begin{enumerate}
  \item Swapping any two rows or any two columns.
  \item Adding integer multiples of one row/column to another row/column.
  \item Multiplying any row/column by $\pm 1$.
\end{enumerate}
Given a square integer matrix $M$, the Smith normal form (SNF) of $M$ is the unique equivalent diagonal matrix $\diag(f_1,f_2,\dots,f_r,0,\dots,0)$ whose non-zero entries are non-negative and satisfy $f_i$ divides $f_{i+1}$, and $r$ is the rank of $M$. 
The diagonal elements of the SNF are known as \emph{invariant factors} or \emph{elementary divisors} of $M$.
In fact, Kannan and Bachem found in \cite{KB} polynomial algorithms for computing the Smith normal form of an integer matrix.
An alternative way of obtaining the Smith normal form is by the \emph{elementary divisors theorem}.
\begin{proposition}
    Let $\Delta_i(M)$ denote the {\it greatest common divisor} of the $i$-minors of the integer matrix $M$, then its $i$-{\it th} invariant factor, $f_i$, is equal to $\Delta_i(M)/ \Delta_{i-1}(M)$, where $\Delta_0(M)=1$.
\end{proposition}

The following observation will give us the relation between the Smith normal form of the distance matrix and the distance ideals.
\begin{proposition}\cite{at}\label{prop:eval1}
Let ${\bf d}\in \mathbb{Z}^{V(G)}$.
If $f_1(M)\mid\cdots\mid f_{r}(M)$ are the invariant factors of the matrix $M:=D_X(G)|_{X={\bf d}}$, then
\[
I_i(G)|_{X={\bf d}}=\langle\Delta_i(M)\rangle=\left\langle \prod_{j=1}^{i} f_j(M) \right\rangle\text{ for all }1\leq i\leq r.
\]
\end{proposition}
Thus to recover $\Delta_i(D(G))$ and the invariant factors of the distance matrix from the distance ideals or the univariate distance ideals, we only need to evaluate $\{I_i(G)\}_{i=1}^n$ at $X={\bf 0}$ or $\{U_i(G)\}_{i=1}^n$ at $t=0$.
Similarly, we can also recover the Smith normal form of the matrices $D^L$, $D^Q$, $D^{\deg}$ or $D^{\deg}_+$.

An ideal is said to be {\it unit} or {\it trivial} if it is equal to $\langle1\rangle$. 
Let $\Phi(G)$ denote the maximum integer $i$ for which the ideal $I_i(G)$ is trivial.
Note that every digraph with at least one non-multiple arc has at least one trivial distance ideals.
On the other hand, let $\phi(M(G))$ denote the number of invariant factors of the matrix $M$ of the digraph $G$ equal to 1, where $M$ is one of the matrices $D$, $D^L$, $D^Q$, $D^{\deg}$ or $D^{\deg}_+$.
Therefore, it follows by Proposition~\ref{prop:eval1} that if the distance ideal $I_i(G)$ is trivial, then $\Delta_i(M)$ and $f_i$ are equal to $1$.
Equivalently, if $\Delta_i(M)$ and $f_i$ are not equal to $1$, then the distance ideal $I_i(G)$ is not trivial.
In this way, distance ideals can be regarded as a generalization of the invariant factors of the matrices $D$, $D^L$, $D^Q$, $D^{\deg}$ or $D^{\deg}_+$.

\begin{corollary}\cite{at}\label{coro:eval1}
    Suppose $M$ is one of the matrices $D$, $D^L$, $D^Q$, $D^{\deg}$ or $D^{\deg}_+$.
    For any digraph $G$, $\Phi(G)\leq \phi(M(G))$.
    And, for any positive integer $k$, the family of digraphs with $\Phi(G)\leq k$ contains the family of digraphs with $\phi(M(G))\leq k$.
\end{corollary}

\section{Strong digraphs with one trivial distance ideal}


Despite that distance ideals are not induce monotone, in \cite{at} there were obtained classifications of the graphs in $\Gamma_1^{\mathbb Z}$ and $\Gamma_1^{\mathbb R}$, which have exactly one trivial distance ideal whose base ring is $\mathbb{Z}$ and $\mathbb{R}$, respectively, in terms of forbidden induced subgraphs.
We will recall them in order to obtain a classification of the digraphs with one trivial distance ideal over $\mathbb{Z}[X]$, that is, digraphs in $\Gamma_1^{\mathbb Z}$.
As previously stated we will only consider strong digraphs

We say that a graph is {\it forbidden} for graphs with at most $k$ trivial distance ideals if any graph containing it as induced subgraph has the $(k+1)$-th distance ideal trivial.
In addition, we say that a forbidden graph $G$ for graphs with at most $k$ trivial distance ideals is {\it minimal} if there is no induced subgraph in $G$ with $k+1$ trivial distance ideals and the $(k+1)$-th distance ideal is trivial for any graph containing $G$ as an induced subgraph.
Next classifications were obtained in \cite{at}, and it is worth mentioning that in \cite{alfaro2} it was found an infinite number of forbidden graphs for the graphs with at most two trivial distance ideals over $\mathbb Z[X]$.

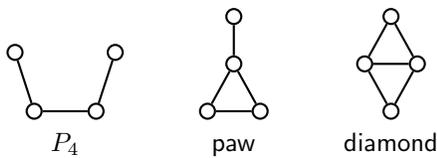
\begin{figure}[h!]
\begin{center}
\begin{tabular}{c@{\extracolsep{10mm}}c@{\extracolsep{10mm}}c@{\extracolsep{10mm}}c@{\extracolsep{10mm}}c}
	\begin{tikzpicture}[scale=.7]
	\tikzstyle{every node}=[minimum width=0pt, inner sep=2pt, circle]
	\draw (126+36:1) node (v1) [draw] {};
	\draw (198+36:1) node (v2) [draw] {};
	\draw (270+36:1) node (v3) [draw] {};
	\draw (342+36:1) node (v4) [draw] {};
	\draw (v1) -- (v2);
	\draw (v2) -- (v3);
	\draw (v4) -- (v3);
	\end{tikzpicture}
&
	\begin{tikzpicture}[scale=.7]
	\tikzstyle{every node}=[minimum width=0pt, inner sep=2pt, circle]
	\draw (-.5,-.9) node (v1) [draw] {};
	\draw (.5,-.9) node (v2) [draw] {};
	\draw (0,0) node (v3) [draw] {};
	\draw (0,.9) node (v4) [draw] {};
	\draw (v1) -- (v2);
	\draw (v1) -- (v3);
	\draw (v2) -- (v3);
	\draw (v3) -- (v4);
	\end{tikzpicture}
&
	\begin{tikzpicture}[scale=.7]
	\tikzstyle{every node}=[minimum width=0pt, inner sep=2pt, circle]
	\draw (-.5,0) node (v2) [draw] {};
	\draw (0,-.9) node (v1) [draw] {};
	\draw (.5,0) node (v3) [draw] {};
	\draw (0,.9) node (v4) [draw] {};
	\draw (v1) -- (v2);
	\draw (v1) -- (v3);
	\draw (v2) -- (v3);
	\draw (v2) -- (v4);
	\draw (v3) -- (v4);
	\end{tikzpicture}
\\
$P_4$
&
$\sf{paw}$
&
$\sf{diamond}$
\end{tabular}
\end{center}
\caption{The graphs $P_4$, $\sf{paw}$ and $\sf{diamond}$.}
\label{fig:forbiddendistance1}
\end{figure}

\begin{proposition}\cite{at}\label{teo:classification}
For $G$ a simple connected graph, the following are equivalent:
\begin{enumerate}
\item $G$ has only 1 trivial distance ideal over $\mathbb{Z}[X]$.
\item $G$ is $\{P_4,\sf{paw},\sf{diamond}\}$-free.
\item $G$ is an induced subgraph of $K_{m,n}$ or $K_{n}$.
\end{enumerate}
\end{proposition}

\begin{proposition}\cite{at}\label{teo:classification2}
For $G$ a simple connected graph, the following are equivalent:
\begin{enumerate}
\item $G$ has only 1 trivial distance ideal over $\mathbb{R}[X]$.
\item $G$ is $\{P_4,\sf{paw},\sf{diamond}, C_4\}$-free.
\item $G$ is an induced subgraph of $K_{1,n}$ or $K_{n}$.
\end{enumerate}
\end{proposition}

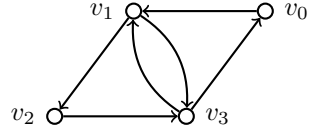
\begin{figure}[ht]
    \centering
    \begin{tikzpicture}[scale=0.7,thick]
	\tikzstyle{every node}=[minimum width=0pt, inner sep=2pt, circle]
        \draw (2,1) node[draw,label={0:$v_0$}] (0) {};
        \draw (-0.5,1) node[draw,label={180:$v_1$}] (1) {};
        \draw (-2,-1) node[draw,label={180:$v_2$}] (2) {};
        \draw (0.5,-1) node[draw,label={0:$v_3$}] (3) {};
        \draw (0) edge[->] (1);
        \draw (1) edge[->] (2);
        \draw (2) edge[->] (3);
        \draw (3) edge[->] (0);
        \draw (1) edge[<-,bend right] (3);
        \draw (3) edge[<-,bend right] (1);
    \end{tikzpicture}
    \caption{A strong digraph whose second distance ideal over $\mathbb Z$ is not trivial.} 
    \label{fig:digraphwithonetrivialideal}
\end{figure}

These classifications can be used in the context of digraphs, that is regarding $K_{m,n}$ or $K_{n}$ as directed graphs, to conclude that these digraphs are included in the family $\Gamma_1^\mathbb{Z}$ of digraphs with one trivial distance ideal over $\mathbb{Z}$.
However, other digraphs in this family do not look like $K_{m,n}$ or $K_{n}$.
For example, consider the digraph shown in Figure~\ref{fig:digraphwithonetrivialideal}, a Gröbner base for the second distance ideal over $\mathbb Z[X]$ of this digraph is $\langle x_0 + 2, x_1 + 1, x_2 + 2, x_3 + 1, 3\rangle$, from which follows that this digraph belongs to $\Gamma_1^{\mathbb Z}$.

For undirected graphs, the forbidden graph $P_4$ bounds the diameter of the graphs in these families, but it does not work for digraphs.
However, we have the following result.

\begin{lemma}\label{lemma:diameterbound}
    Let $G$ be a strong digraph. 
    If $\diam(G)\geq 3$ then $I_2(G)= \langle 1 \rangle.$
\end{lemma}
\begin{proof}
    Assume $u,v\in V(G)$ such that $\dist(u,v)=3$. 
    Then there exist vertices $w,z$ such that $uw, wz, zv \in A(D)$ and $uz, wv, vz \notin A(D)$. 
    This means that the submatrix of the distance matrix generated by the rows of $w$ and $u$, and the columns of $v$ and $z$  satisfies that

    \[
    \det
    \begin{bmatrix}
        \dist(w,v) & \dist(w,z) \\
        \dist(u,v) & \dist(u,z) \\
    \end{bmatrix}
    =
    \det
    \begin{bmatrix}
        2 & 1 \\
        3 & 2 \\
    \end{bmatrix}
    = 1.
    \]
    From which the result follows.
\end{proof}

The previous lemma implies that strong digraphs in $\Gamma_1^{\mathbb Z}$ have diameter at most two. 
This also can be regarded as a forbidden combinatorial structure within the digraph.
In Section~\ref{sec:circuits}, we will see that the circuits of length at least 4 are minimal forbidden subdigraphs for the graphs in $\Gamma_{\leq1}^\mathbb{Z}$.
However, we will introduce a different combinatorial structure, which we will call {\it pattern}, in order to obtain a classification of the strong digraphs in $\Gamma_{\leq1}^\mathbb{Z}$.
To settle this point properly let us introduce this concept. 

\begin{definition}
Let $D=(V,A)$ a strong digraph, a {\it pattern} ${\mathcal P}=(U,B,C)$ in $D$ is the subset of vertices $U\subseteq V$ together with two disjoint sets of arcs, $B$ and $C$, whose endpoints are in $U$, where arcs in $B$ are arcs in $A$ while arcs in $C$ must not belong to $A$.    
\end{definition}



\begin{figure}
    \centering
    \begin{tabular}{c@{\extracolsep{2cm}}cc}
    \begin{tikzpicture}[scale=1.7]
        \node[label={180:$u$}] (z) at (0,0){};
        \node[label={0:$w$}] (w) at (1,0){};
        \node[label={0:$z$}] (y) at (1,1){};
        \node[label={180:$v$}] (x) at (0,1){};
        \draw[->] (z) to (w);
        \draw[->] (w) to (y);
        \draw[->] (y) to (x);
        \draw[->, dotted] (z) to (y);
        \draw[->, dotted] (z) to (x);
        \draw[->, dotted] (w) to (x);
    \end{tikzpicture}
    &
    \begin{tikzpicture}[scale=1.7]
        \node[label={180:$u$}] (z) at (0,0){};
        \node[label={0:$v$}] (w) at (1,0){};
        \node[label={0:$w$}] (y) at (1,1){};
        \node[label={180:$z$}] (x) at (0,1){};
        \draw[->] (z) to (w);
        \draw[->] (w) to (y);
        \draw[->, dotted] (z) to (y);
        \draw[->, dotted] (z) to (x);
        \draw[->, dotted] (w) to (x);
    \end{tikzpicture}
        &
    \begin{tikzpicture}[scale=1,thick]
    \tikzstyle{every node}=[minimum width=0pt, inner sep=2pt, circle]
        \draw (0:1) node[draw,label={0:$u$}] (0) { };
        \draw (120:1) node[draw,label={180:$v$}] (1) { };
        \draw (240:1) node[draw,label={180:$w$}] (2) { };
        \draw  (0) edge[->,bend right=10] (1);
        \draw  (0) edge[->,bend right=10] (2);
        \draw  (1) edge[->,dotted,bend right=10] (0);
        \draw  (1) edge[->,bend right=10] (2);
        \draw  (2) edge[->,bend right=10] (0);
        \draw  (2) edge[->,bend right=10] (1);
    \end{tikzpicture}
    \\
    $\mathcal{P}_4$ & $\mathcal{F}_1$ & $\mathcal{F}_2$
    \end{tabular}

    \caption{The forbidden patterns $\mathcal{P}_4$, $\mathcal{F}_1$ and $\mathcal{F}_2$, in which the dotted arrow means that $vu\notin A(G)$.}
    \label{fig:forbiddenstructure3vertices}
\end{figure}
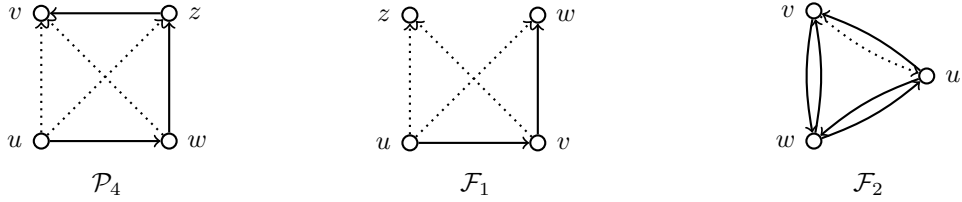

A pattern allows us to describe forbidden structures within digraphs forbidden for digraphs with one trivial distance ideal.
For example, consider the patterns described in Figure~\ref{fig:forbiddenstructure3vertices}, where arcs in $B$ are solid while arcs in $C$ are dotted; the arcs not described by a pattern are not important.
Note that the proof of Lemma~\ref{lemma:diameterbound} also implies that if a digraph whose diameter is at least 3, then it contains a pattern $\mathcal{P}_4$, at the left-hand side in Figure~\ref{fig:forbiddenstructure3vertices}.
This does not imply that every pattern $\mathcal{P}_4$ within the digraph implies that the second distance ideal is trivial, this is because the distance from $u$ to $v$ in the digraph might be 2. 
However, Lemma~\ref{lemma:diameterbound} allows us to assume that strong digraphs whose second distance ideal is not trivial have diameter 2.
With this in mind, let us prove that the patterns $\mathcal{F}_1$ and $\mathcal{F}_2$ of Figure~\ref{fig:forbiddenstructure3vertices} are not contained in any strong digraph with at most one trivial distance ideal.

\begin{lemma}
    If $G$ is a strong digraph with diameter at most 2 that contains the pattern $\mathcal{F}_1$ of Figure~\ref{fig:forbiddenstructure3vertices}, then $I_2(G)=\langle 1 \rangle$.
\end{lemma}
\begin{proof}
    If $G$ contains pattern $\mathcal{F}_1$, then $D_X(G)$ contains the 2-minors:
    \[
    \det
    \begin{bmatrix}
        \dist(u,v) & \dist(u,w)\\
        \dist(v,v) & \dist(v,w)\\
    \end{bmatrix}
    =
    \begin{bmatrix}
        1 & 2\\
        x_v & 1\\
    \end{bmatrix}
    = 1-2x_v
    \]
    and
    \[
    \det
    \begin{bmatrix}
        \dist(u,v) & \dist(u,z)\\
        \dist(v,v) & \dist(v,z)\\
    \end{bmatrix}
    =
    \begin{bmatrix}
        1 & 2\\
        x_v & 2\\
    \end{bmatrix}
    = 2-2x_v
    \]
    From which follows that $1\in I_2(G)$.
\end{proof}

On the other hand, note that any digraph of diameter at least 3 contains a pattern $\mathcal{F}_1$, this is because pattern $\mathcal{P}$ contains a pattern isomorphic to pattern $\mathcal{F}_1$. Therefore, if $\mathcal{F}_1$ is a forbidden pattern for $G$, then $\diam (G) \leq 2$.

\begin{lemma} \label{P}
    If $G$ is a strong digraph with diameter at most 2 that contains the pattern $\mathcal{F}_2$ of Figure~\ref{fig:forbiddenstructure3vertices}, then $I_2(G)=\langle 1 \rangle$. 
\end{lemma}
\begin{proof}
    Note the submatrix 
    $D_X(P)=
    \begin{bmatrix}
        x_u & 1   & 1   \\
        2   & x_v & 1   \\
        1   & 1   & x_w \\
    \end{bmatrix}
    $ of $D_X(G)$ contains the 2-minors:
    $\det
    \begin{bmatrix}
        x_u & 1 \\
        2   & 1 \\
    \end{bmatrix}=x_u-2$
    and
    $\det
    \begin{bmatrix}
        x_u & 1 \\
        1   & 1 \\
    \end{bmatrix}=x_u-1$.
    Therefore $1\in I_2(G)$, and the result turns out.
\end{proof}
Now let us consider patterns having four vertices.
\begin{lemma} \label{34}
    If $G$ is a strong digraph with diameter at most 2 that contains the pattern $\mathcal{F}_3$ of Figure~\ref{3p}, then $I_2(G)=\langle 1 \rangle$.
\end{lemma}
\begin{figure}[h!]
\centering
\begin{tabular}{c@{\extracolsep{2cm}}ccc}
    \begin{tikzpicture}[scale=1.5]
        \node[label={180:$w$}] (z) at (0,0){};
        \node[label={0:$u$}] (w) at (1,0){};
        \node[label={0:$z$}] (y) at (1,1){};
        \node[label={180:$v$}] (x) at (0,1){};
        \draw[->] (w) to (y);
        \draw[->] (w) to (x);
        \draw[->] (z) to (x);
        \draw[->, dotted] (z) to (y);
    \end{tikzpicture}
    &
    \begin{tikzpicture}[scale=1]
        \node[label={0:$w$}] (w) at (-30:1){};
        \node[label={0:$v$}] (v) at (90:1){};
        \node[label={180:$u$}] (u) at (210:1){};
        \draw[->,bend right=15] (u) to (v);
        \draw[->,bend right=40] (v) to (u);
        \draw[->,bend right=15] (v) to (w);
        \draw[->,bend right=40,dotted] (w) to (v);
        \draw[->,bend right=40,dotted] (u) to (w);
        \draw[->,bend right=15,dotted] (w) to (u);
    \end{tikzpicture}
    &
    \begin{tikzpicture}[scale=1]
        \node[label={0:$w$}] (w) at (-30:1){};
        \node[label={0:$v$}] (v) at (90:1){};
        \node[label={180:$u$}] (u) at (210:1){};
        \draw[->,bend right=15] (u) to (v);
        \draw[->,bend right=40] (v) to (u);
        \draw[->,bend right=15,dotted] (v) to (w);
        \draw[->,bend right=40] (w) to (v);
        \draw[->,bend right=40,dotted] (u) to (w);
        \draw[->,bend right=15,dotted] (w) to (u);
    \end{tikzpicture}
    \\
    $\mathcal{F}_3$ & $\mathcal{F}_4$ & $\mathcal{F}_5$
\end{tabular}
\caption{Three forbidden patterns.}\label{3p}
\end{figure}
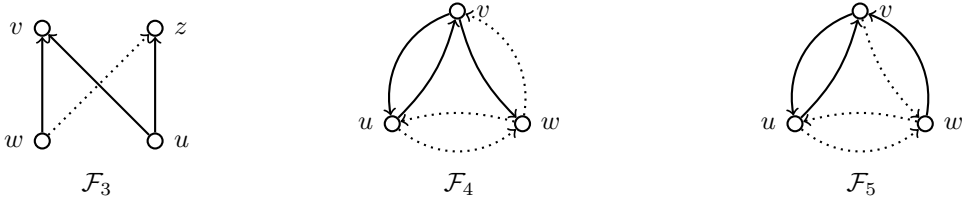
\begin{proof}
    Suppose $G$ contains the patterm $\mathcal{F}_3$, then $D_X(G)$ contains the submatrix 
    \[
    \begin{bmatrix}
        \dist(u,v) & \dist(u,z) \\
        \dist(w,v) & \dist(w,z)
    \end{bmatrix}
    =
    \begin{bmatrix}
        1 & 1 \\
        1 & 2
    \end{bmatrix},
    \]
    which has determinant equal to $1$.
\end{proof}

\begin{lemma}\label{forb}
    If $G$ is a strong digraph with diameter at most 2 that contains the patterns $\mathcal{F}_4$ and $\mathcal{F}_5$ of Figure~\ref{3p}, then $I_2(G)=\langle 1 \rangle$.
\end{lemma}
\begin{proof}
    Note if $G$ contains pattern $\mathcal{F}_4$, then $D_X(G)$ contains the 2-minors:
    \[
    \det
    \begin{bmatrix}
        \dist(u,v) & \dist(u,w)\\
        \dist(v,v) & \dist(v,w)\\
    \end{bmatrix}
    =
    \begin{bmatrix}
        1 & 2\\
        x_v & 1\\
    \end{bmatrix}
    = 1-2x_v
    \]
    and
    \[
    \det
    \begin{bmatrix}
        \dist(v,u) & \dist(v,v)\\
        \dist(w,u) & \dist(w,v)\\
    \end{bmatrix}
    =
    \begin{bmatrix}
        1 & x_v\\
        2 & 2\\
    \end{bmatrix}
    = 2-2x_v.
    \]
    From which follows that $1\in I_2(G)$.
    A similar argument implies that a strong digraph with diameter at most 2 that contains pattern $\mathcal{F}_5$ has $I_2(G)=\langle1\rangle$.
\end{proof}

Now, we are ready to state the main theorem of this section.

\begin{theorem}\label{equi}
    Let $G$ be a strong digraph. 
    The following statements are equivalent:
    \begin{enumerate}
        \item $G$ has exactly one trivial distance ideal over $\mathbb Z[X]$,
        \item Patterns $\mathcal{F}_1$, $\mathcal{F}_2$, $\mathcal{F}_3$, $\mathcal{F}_4$ and $\mathcal{F}_5$ are forbidden for $G$. 
        \item $G$ is one of the strong digraphs described by the following diagrams:
        \begin{center}
        \begin{tabular}{c@{\extracolsep{2cm}}c}
        \begin{tikzpicture}
            \node[fill=black!20] (z) at (0,0){$K_a$};
            \node (y) at (1,1){$T_d$};
            \node (x) at (1,-1){$T_b$};
            \node[fill=black!20] (w) at (2,0){$K_c$};
            \draw[->] (z)--(x);
            \draw[->] (y)--(z);
            \draw[->] (x)--(w);
            \draw[->] (w)--(y);
            \draw[->, bend left] (x) to (y);
            \draw[->, bend left] (y) to (x);
        \end{tikzpicture}
        &
        \begin{tikzpicture}
            \clip (-1.4,-1.4) rectangle (1.4,1.4);
            \node (a) at (0:1){};
            \node (b) at (120:1){};
            \node (c) at (240:1){};
            \draw[->] (a)--(b);
            \draw[->] (b)--(c);
            \draw[->] (c)--(a);
        \end{tikzpicture}\\
        $\Lambda(a,b,c,d)$ & $\overrightarrow{C_3}$
        \end{tabular}
        \end{center}
    \end{enumerate}
    where $T_p$ denote an independent set of order $p\geq 0$, and $K_q$ denote a complete digraph of order $q\geq0$, and an arrow between two sets, say from $A$ to $B$, means that there is an arrow from each vertex in $A$ to each vertex of $B$.
\end{theorem}
\begin{proof}
    $(1) \Rightarrow (2):$ Lemma 7 implies that the diameter of $G$ is at most 2. 
    If $G$ has exactly one trivial distance ideal over $\mathbb{Z}$, then, by Lemmas \ref{P}, \ref{34} and \ref{forb}, the patterns $\mathcal{F}_1, \dots, \mathcal{F}_5$ are forbidden for $G$.
    
    $(2) \Rightarrow (3):$ We will proceed by induction over $|V(G)|$. 
    In Figure~\ref{fig:lambdaupto4}, we can see the strong digraphs with up to 4 vertices, except by $\overrightarrow{C_3}$, containing none of the forbidden patterns. 
    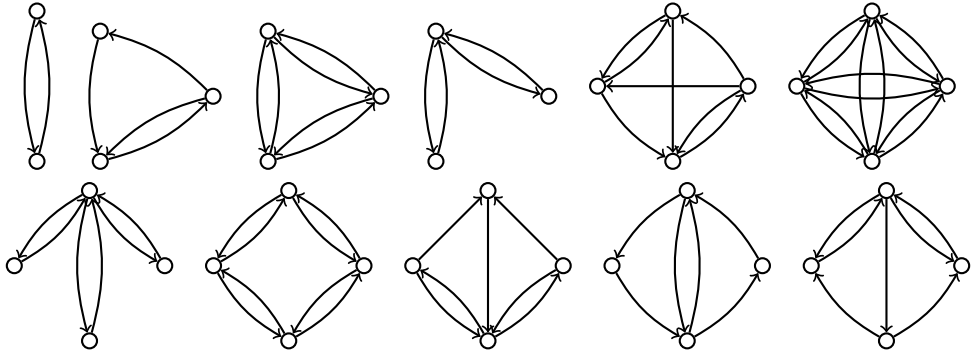
\begin{figure}[h]
        \centering
        \begin{tabular}{cccccc}
            \begin{tikzpicture}[scale=1,thick]
            \tikzstyle{every node}=[minimum width=0pt, inner sep=2pt, circle]
                \draw (90:1) node[draw] (0) {}; 
                \draw (270:1) node[draw] (1) {}; 
                \draw[->,bend right=15]  (0) edge (1);
                \draw[->,bend right=15]  (1) edge (0);
            \end{tikzpicture}
            &
            \begin{tikzpicture}[scale=1,thick]
            \tikzstyle{every node}=[minimum width=0pt, inner sep=2pt, circle]
                \draw (0:1) node[draw] (0) {}; 
                \draw (120:1) node[draw] (1) {}; 
                \draw (240:1) node[draw] (2) {}; 
                \draw[->,bend right=15]  (0) edge (1);
                \draw[->,bend right=15]  (0) edge (2);
                \draw[->,bend right=15]  (1) edge (2);
                \draw[->,bend right=15]  (2) edge (0);
            \end{tikzpicture}
            &
            \begin{tikzpicture}[scale=1,thick]
            \tikzstyle{every node}=[minimum width=0pt, inner sep=2pt, circle]
                \draw (0:1) node[draw] (0) {}; 
                \draw (120:1) node[draw] (1) {}; 
                \draw (240:1) node[draw] (2) {}; 
                \draw[->,bend right=15]  (0) edge (1);
                \draw[->,bend right=15]  (0) edge (2);
                \draw[->,bend right=15]  (1) edge (0);
                \draw[->,bend right=15]  (1) edge (2);
                \draw[->,bend right=15]  (2) edge (0);
                \draw[->,bend right=15]  (2) edge (1);
            \end{tikzpicture}
            &
            \begin{tikzpicture}[scale=1,thick]
            \tikzstyle{every node}=[minimum width=0pt, inner sep=2pt, circle]
                \draw (0:1) node[draw] (0) {}; 
                \draw (120:1) node[draw] (1) {}; 
                \draw (240:1) node[draw] (2) {}; 
                \draw[->,bend right=15]  (0) edge (1);
                \draw[->,bend right=15]  (1) edge (0);
                \draw[->,bend right=15]  (1) edge (2);
                \draw[->,bend right=15]  (2) edge (1);
            \end{tikzpicture}
            &
            \begin{tikzpicture}[scale=1,thick]
            \tikzstyle{every node}=[minimum width=0pt, inner sep=2pt, circle]
                \draw (0:1) node[draw] (0) {}; 
                \draw (90:1) node[draw] (1) {}; 
                \draw (180:1) node[draw] (2) {}; 
                \draw (270:1) node[draw] (3) {}; 
                \draw[->,bend right=15]  (0) edge (1);
                \draw[->]  (0) edge (2);
                \draw[->,bend right=15]  (0) edge (3);
                \draw[->,bend right=15]  (1) edge (2);
                \draw[->]  (1) edge (3);
                \draw[->,bend right=15]  (2) edge (1);
                \draw[->,bend right=15]  (2) edge (3);
                \draw[->,bend right=15]  (3) edge (0);
            \end{tikzpicture}
            &
            \begin{tikzpicture}[scale=1,thick]
            \tikzstyle{every node}=[minimum width=0pt, inner sep=2pt, circle]
                \draw (0:1) node[draw] (0) {}; 
                \draw (90:1) node[draw] (1) {}; 
                \draw (180:1) node[draw] (2) {}; 
                \draw (270:1) node[draw] (3) {}; 
                \draw[->,bend right=15]  (0) edge (1);
                \draw[->,bend right=15]  (0) edge (2);
                \draw[->,bend right=15]  (0) edge (3);
                \draw[->,bend right=15]  (1) edge (0);
                \draw[->,bend right=15]  (1) edge (2);
                \draw[->,bend right=15]  (1) edge (3);
                \draw[->,bend right=15]  (2) edge (0);
                \draw[->,bend right=15]  (2) edge (1);
                \draw[->,bend right=15]  (2) edge (3);
                \draw[->,bend right=15]  (3) edge (0);
                \draw[->,bend right=15]  (3) edge (1);
                \draw[->,bend right=15]  (3) edge (2);
            \end{tikzpicture}
        \end{tabular}
        \begin{tabular}{ccccc}
            \begin{tikzpicture}[scale=1,thick]
            \tikzstyle{every node}=[minimum width=0pt, inner sep=2pt, circle]
                \draw (0:1) node[draw] (0) {}; 
                \draw (90:1) node[draw] (1) {}; 
                \draw (180:1) node[draw] (2) {}; 
                \draw (270:1) node[draw] (3) {}; 
                \draw[->,bend right=15]  (0) edge (1);
                \draw[->,bend right=15]  (1) edge (0);
                \draw[->,bend right=15]  (1) edge (2);
                \draw[->,bend right=15]  (1) edge (3);
                \draw[->,bend right=15]  (2) edge (1);
                \draw[->,bend right=15]  (3) edge (1);
            \end{tikzpicture}
            &
            \begin{tikzpicture}[scale=1,thick]
            \tikzstyle{every node}=[minimum width=0pt, inner sep=2pt, circle]
                \draw (0:1) node[draw] (0) {}; 
                \draw (90:1) node[draw] (1) {}; 
                \draw (180:1) node[draw] (2) {}; 
                \draw (270:1) node[draw] (3) {}; 
                \draw[->,bend right=15]  (0) edge (1);
                \draw[->,bend right=15]  (0) edge (3);
                \draw[->,bend right=15]  (1) edge (0);
                \draw[->,bend right=15]  (1) edge (2);
                \draw[->,bend right=15]  (2) edge (1);
                \draw[->,bend right=15]  (2) edge (3);
                \draw[->,bend right=15]  (3) edge (0);
                \draw[->,bend right=15]  (3) edge (2);
            \end{tikzpicture}
            &
            \begin{tikzpicture}[scale=1,thick]
            \tikzstyle{every node}=[minimum width=0pt, inner sep=2pt, circle]
                \draw (0:1) node[draw] (0) {}; 
                \draw (90:1) node[draw] (1) {}; 
                \draw (180:1) node[draw] (2) {}; 
                \draw (270:1) node[draw] (3) {}; 
                \draw[->]  (0) edge (1);
                \draw[->,bend right=15]  (0) edge (3);
                \draw[->]  (1) edge (3);
                \draw[->]  (2) edge (1);
                \draw[->,bend right=15]  (2) edge (3);
                \draw[->,bend right=15]  (3) edge (0);
                \draw[->,bend right=15]  (3) edge (2);
            \end{tikzpicture}
            &
            \begin{tikzpicture}[scale=1,thick]
            \tikzstyle{every node}=[minimum width=0pt, inner sep=2pt, circle]
                \draw (0:1) node[draw] (0) {}; 
                \draw (90:1) node[draw] (1) {}; 
                \draw (180:1) node[draw] (2) {}; 
                \draw (270:1) node[draw] (3) {}; 
                \draw[->,bend right=15]  (0) edge (1);
                \draw[->,bend right=15]  (1) edge (2);
                \draw[->,bend right=15]  (1) edge (3);
                \draw[->,bend right=15]  (2) edge (3);
                \draw[->,bend right=15]  (3) edge (0);
                \draw[->,bend right=15]  (3) edge (1);
            \end{tikzpicture}
            &
            \begin{tikzpicture}[scale=1,thick]
            \tikzstyle{every node}=[minimum width=0pt, inner sep=2pt, circle]
                \draw (0:1) node[draw] (0) {}; 
                \draw (90:1) node[draw] (1) {}; 
                \draw (180:1) node[draw] (3) {}; 
                \draw (270:1) node[draw] (2) {}; 
                \draw[->,bend right=15]  (0) edge (1);
                \draw[->,bend right=15]  (1) edge (0);
                \draw[->]  (1) edge (2);
                \draw[->,bend right=15]  (1) edge (3);
                \draw[->,bend right=15]  (2) edge (0);
                \draw[->,bend left=15]  (2) edge (3);
                \draw[->,bend right=15]  (3) edge (1);
            \end{tikzpicture}
        \end{tabular}
        \caption{The strong digraphs, except by $\overrightarrow{C_3}$ with up to 4 vertices containing none of the patterns $\mathcal{F}_1$, $\mathcal{F}_2$, $\mathcal{F}_3$, $\mathcal{F}_4$ and $\mathcal{F}_5$}
        \label{fig:lambdaupto4}
    \end{figure}
    It is easy to see that a complete $n$-partite digraph contains a forbidden pattern $\mathcal{F}_3$ for $n \geq 3.$
    Note that complete and complete bipartite digraphs do not contain the forbidden patterns and they can be described as $K_a=\Lambda(a,0,0,0), K_{b,d}=\Lambda(0,b,0,d)$.
    Thus, we will assume $G$ is not a complete or a complete bipartite digraph, this means that there exist vertices $v$ and $u$ such that that $(vu)\in A(G)$ and $(uv)\notin A(G)$.

    
    By induction hypothesis we know that $G\setminus v$ is of the form $\Lambda(a,b,c,d)$.

    Assume first that $ u \in K_a,$ let $ t\in T_d, i\in T_b$ and $ k \in K_c$. 
    Notice that we have at least one of the arcs $(vk)$ and $(vi)$ since otherwise we obtain the forbidden pattern $\mathcal{P}_4$.
    If $(vk)\in A(G)$, the vertices $v,k,u,t$ form a forbidden pattern $\mathcal{F}_3$. 
    This means that $(vi)\in A(G)$. 
    Since $G$ is strongly connected of diameter two and $(uv)\notin A(G)$, we have that $(iv)$ for some $ i \in T_b$. 
    This means that there are no arcs between $v$ and $T_d$ since otherwise we obtain the forbidden pattern $\mathcal{F}_2$. 
    This also means that we must have every arc from $T_b$ to $v$, since otherwise we would obtain the forbidden pattern $\mathcal{F}_5$. 
    Finally, since there are no arcs from $T_d$ to $v$, we must have $(kv)\in A(G)$ since otherwise $d(k,v)>2$.
 
    If there exists another vertex $w \in K_a$ then it is not possible that there are no arcs between $w$ and $v$, since otherwise $\mathcal{F}_5$ appears, and it is not possible that both arcs $(wv)$ and $(vw)$ exists since ${\mathcal F}_2$ would appear. 
    Thus either $(vw)$ or $(wv)\in A(G)$. 
    If $(wv)\in A(G)$ we obtain the forbidden pattern $\mathcal{F}_3$  with another  vertex in  vertex in $T_b$.

    This means that $(tv)\in A(G)$ and thus $G=\Lambda(a,b,c,d+1)$. 
    If $u \in K_c, $ analogously we obtain that $G=\Lambda(a,b+1,c,d)$.

    Assume that $u\in T_b$ and let $t\in T_d$. 
    If there are no arcs between $t$ and $v$, we obtain the forbidden pattern $\mathcal{F}_5$. 
    If $(vt)\in A(G)$ we obtain the forbidden pattern $\mathcal{F}_3$ with a vertex of $K_c$. 
    This means that $(tv)\in A(G)$ for every $t \in T_d$.
 

    Let $k \in K_{a}$ and notice that if $(vk)\notin A(G)$ we obtain the forbidden pattern ${\mathcal F}_3$. 
    Since $kutv$ is a path of length four with $(uv), (kt)\notin A(G)$, we have that $(kv)\in A(G)$ thus $v$ must be adjacent to every vertex in $K_{a}$. 
    Let $w$ be any other vertex in $T_b$. 
    Then $vutw$ is a path of length four such that $(uw), (vt)\notin A(G)$ thus $(vw)\in A(G)$. 
    Any arc from $T_b$ to $v$ leads to the forbidden pattern ${\mathcal F}_2$. 
    Finally notice that any arc between $K_{c}$ and $v$ leads to a forbidden pattern ${\mathcal F}_4$ or ${\mathcal F}_5$ if $a>0$ and leads to the forbidden pattern ${\mathcal F}_{3}$ if $a=0$. This means that $G=\Lambda(a+1,b,c,d)$ and if $u \in T_d$ we obtain analogously that $G=\Lambda(a,b,c+1,d)$.

    $(3)\Rightarrow (1):$ 
     The matrix $D_X(\overrightarrow{C}_{3})$ is
     \[
    \begin{bmatrix}
       x_1 & 1 & 2  \\ 
       2 & x_2 & 1  \\ 
       1 & 2 & x_3 \\
    \end{bmatrix}.
    \]
    Since there are entries equal to $1$, then the first distance ideal is trivial.
    The $2$-minors of $D_X(\overrightarrow{C}_{3})$ are $x_1 x_2 - 2, x_1 - 4, -2 x_2 + 1, 2 x_1 - 1, x_1 x_3 - 2, x_3 - 4, -x_2 + 4, 2 x_3 - 1, x_2 x_3 - 2$.
    A Gröbner basis of the second distance ideal is $\langle x_1 + 3, x_2 + 3, x_3 + 3, 7\rangle$.
    
    On the other hand, the matrix $D_X$ of $\Lambda(a,b,c,d)$ is 
    \[
    \diag(x_1,\dots,x_a,x_{a+1},\dots,x_{a+c},y_{1},\dots,y_{b},y_{b+1},\dots,y_{b+d})+
    \begin{bmatrix}
        {\sf J}_a - {\sf I}_a & 2{\sf J}_{a,c} & {\sf J}_{a,b} & 2{\sf J}_{a,d}\\
        2{\sf J}_{c,a} & {\sf J}_c - {\sf I}_c & 2{\sf J}_{c,b} & {\sf J}_{c,d}\\
        2{\sf J}_{b,a} & {\sf J}_{b,c} & 2{\sf J}_b - 2{\sf I}_b & {\sf J}_{b,d}\\
        {\sf J}_{d,a} & 2{\sf J}_{d,c} & {\sf J}_{d,b} & 2{\sf J}_d - 2{\sf I}_d
    \end{bmatrix},
    \]
    where $\sf I$ and $\sf J$ denote identity matrix and the all-ones matrix, respectively.
    Then the second distance ideal of $\Lambda(a,b,c,d)$, when it is a strong digraph, is the following:
    \begin{equation}
    I_{2}(\Lambda(a,b,c,d))=
        \begin{cases}
            \langle x_{1}x_{2}-1\rangle & \text{if } a=2,~b=c=d=0\\
            \langle x_{1}x_{2}-1\rangle & \text{if } c=2,~a=b=d=0\\
            \langle \{x_i-1\}_{i=1}^{a}\rangle & \text{if } 3\leq a,~ b=c=d=0\\
            \langle \{x_i-1\}_{i=1}^{c}\rangle & \text{if } 3\leq c,~ a=b=d=0\\
            \langle y_{1}y_{2}-1\rangle & \text{if } a=c=0,~b=d=1\\
            \langle \{y_j-2\}_{j=1}^{b}, 2y_{b+1}-1\rangle & \text{if } a=c=0,~ d=1,~b\geq 2\\
            \langle 2y_{1}-1,\{y_j-2\}_{j=2}^{d+1}\rangle & \text{if } a=c=0,~ b=1,~d\geq 2\\
            \left\langle 3, \{x_i-1\}_{i=1}^{a+c}, \{y_j-2\}_{j=1}^{b+d} \right\rangle & \text{otherwise.}
        \end{cases}
    \end{equation}
    Observe that any linear combination of the polynomials $3, \{x_i-1\}_{i=1}^{a+c}, \{y_j-2\}_{j=1}^{b+d}$, when evaluated at $x_{i}=1$ and $y_{j}=2$ for all $i=1, \dots, a+c$ and $j=1, \dots, b+d$, results in a multiple of 3. 
    Consequently, the ideal $\left\langle 3, \{x_i-1\}_{i=1}^{a+c}, \{y_j-2\}_{j=1}^{b+d} \right\rangle$ is not trivial.
    Similarly, the other cases are non-trivial.
\end{proof}

Note that the strong digraphs described in $\Lambda(a,b,c,d)$ contain the complete graphs and the complete bipartite graphs, regarded as digraphs.
This was expected since the characterization of undirected connected graphs with only one trivial distance ideal coincides with such graphs.

In the following we will give an alternative proof that sentence (2) implies sentence (3) of Theorem \ref{equi}.

\begin{proof}[Proof of (2) implies (3) of Theorem \ref{equi}]
We will proceed by induction over the number of vertices. 
In Figure 7, we can see all strong digraphs with at most four vertices, except by $\overrightarrow{C_3}$, such that only the first distance ideal is trivial. 
Each of these graphs is either $\overrightarrow{C_3}$ or can be expressed as $\Lambda(a,b,c,d)$.
Now take $G$ with $n\geq 5$ vertices, thus by induction $G\setminus v$ has the form $G\setminus v = \Lambda (a,b,c,d)$.
First, if $G\setminus v$ is a complete digraph, then $G$ is also a complete digraph by strong connectivity and the forbidden patterns ${\mathcal F}_2$ and ${\mathcal F}_3$. 
Also, by strong connectivity, it is easy to see that if $G\setminus v$ is an edgeless digraph then $G$ is the complete bipartite graph $\overrightarrow{K}_{1,(n-1)}$. Moreover, note that $G\setminus v\neq \Lambda (a,0,c,0)$ with $a,c\geq 1$ since this would form a directed paw and, consequently, would have the pattern ${\mathcal F}_3$. \\
Now we consider the case $(a,b,c,d)=(a,b,0,0)$ which is analogous to the case $(0,0,c,d)$. 
Note that any vertex in $T_b$ has an arc towards $v$ and that $v$ has an arc towards at least one vertex in $K_a$ by strong connectivity. 
Note that since ${\mathcal F}_1$ is forbidden we have that $G$ has diameter at most 2.
If $a\geq 2$, $v$ has an arc towards every vertex in $K_a$ since the diameter is at most two. 
In order to avoid the pattern ${\mathcal F}_2$, every vertex in $K_a$ has an arc towards $v$ or none of them does. 
In the former case, $G$ forms an ${\mathcal F}_3$. 
Therefore, there are no arcs from $K_a$ to $v$. 
Moreover, there is an arc from $v$ to every vertex in $T_b$ otherwise we also have an ${\mathcal F}_3$. 
Thus $G=\Lambda(a,b,0,1)$.
On the other hand, if $a=1$, $v$ has at least one arc towards a vertex in $T_b$, otherwise we have an ${\mathcal F}_1$. 
Thus we have every arc from $v$ to $T_b$ or the graph yields an ${\mathcal F}_5$. In this case, there is no arc from $K_a$ to $v$ or we build an ${\mathcal F}_2$ pattern, concluding that $G=(a,b,0,1)$ as well.\\
If $(a,b,c,d)=(a,0,0,d)$ which is analogous to $(0,b,c,0)$.
Note that $v$ has an arc towards every vertex in $T_d$ and that every vertex in $K_a$ has an arc towards at least one vertex in $K_a$ by strong connectivity and the fact that we have diameter 2.
If $a\geq 2$, there is an arc from $v$ towards every vertex in $K_a$ or there is no such arc to avoid the pattern ${\mathcal F}_2$. 
In the former case, $G$ forms an ${\mathcal F}_3$. 
Therefore, there are no arcs from $v$ to $K_a$. Moreover, there is an arc from $v$ to every vertex in $T_b$ otherwise we also have a pattern ${\mathcal F}_3$. 
Thus $G=\Lambda(a,1,0,d)$.
On the other hand, if $a=1$, $v$ has at least one arc coming from a vertex in $T_d$, otherwise we have the pattern ${\mathcal F}_1$. 
Thus we have every arc from $T_d$ to $v$ or the graph yields an ${\mathcal F}_4$. 
In this case, there is no arc from $v$ to $K_a$ or we build an ${\mathcal F}_2$, thus $G=(a,1,0,d)$.\\
If $(a,b,c,d)=(0,b,0,d)$ with $b,d\geq 1$ and let us assume that $b\geq d$. 
If $v$ is both-ways adjacent to a vertex in $T_d$, then it is either non-adjacent or both-ways adjacent to any vertex in $T_b$ otherwise, we have the pattern ${\mathcal F}_2$. 
If $v$ is adjacent to a vertex in $T_b$ then it is adjacent to every vertex in $T_b$ or we would have a directed paw, but in this case, we have a directed diamond which also implies the pattern ${\mathcal F}_3$. 
Thus $G=\Lambda(0,b+1,0,d)$. 
Similarly, if $v$ is both-ways adjacent to a vertex in $T_b$, we have $G=\Lambda(0,b,0,d+1)$. 
Assuming $v$ is not twice adjacent to any vertex in $G\setminus v$, then by strong connectivity we have $G=\Lambda(1,b,0,d)$ or $G=\Lambda (0,b,1,d)$.\\ 
Every case where we only have one parameter equal to zero can be proven by a recursive argument using the previous cases and considering the forbidden patterns.\\
Finally, if $G\setminus v= \Lambda (a,b,c,d)$ and none of the parameters is zero. 
Assume there is a vertex in $K_a$ with both arcs to $v$. 
If there is another vertex in $K_a$, $v$ needs to be both-ways adjacent to it or we get either an ${\mathcal F}_2$ or an ${\mathcal F}_1$ adding a vertex in $K_c$. 
Then, avoiding ${\mathcal F}_2$, ${\mathcal F}_4$ and ${\mathcal F}_5$, $v$ has to have only one arc from or to $T_b$ ($T_d$). 
Therefore there is an arc from each vertex in $T_d$ to $v$ otherwise we have an ${\mathcal F}_3$. 
Also, there is an arc from $v$ to every vertex of $T_b$. 
Moreover, $v$ is not adjacent to any vertex in $K_c$ to avoid ${\mathcal F}_3$, ${\mathcal F}_4$ and ${\mathcal F}_5$, thus $G=\Lambda (a+1,b,c,d)$. 
On the other hand, if $v$ has only the arc towards a vertex in $K_a$, then that is the case with any other vertex in $K_a$ otherwise an ${\mathcal F}_2$, an ${\mathcal F}_1$ or an ${\mathcal F}_5$ forms. 
Also, $v$ is not adjacent to any vertex of $T_d$, or an ${\mathcal F}_3$ is formed. 
Moreover, both arcs from $v$ and any vertex in $T_b$ are present since otherwise we would have an ${\mathcal F}_1$, an ${\mathcal F}_4$, or an ${\mathcal F}_5$ in $G$. 
Furthermore, every vertex in $K_c$ has only the arc towards $v$ to avoid ${\mathcal F}_4$ and ${\mathcal F}_3$, thus $G=\Lambda (a,b,c,d+1)$. 
Similar arguments hold whether there is a vertex in $K_a$ with only the arc in direction to $v$ ($G=\Lambda (a,b+1,c,d$) or a vertex in $K_a$ non-adjacent to $v$ ($G=\Lambda (a,b,c+1,d)$).
\end{proof}

Next characterizations follow by evaluating the second distance ideals of $\overrightarrow{C}_3$ and $\Lambda(a,b,c,d)$ properly. 
These ideals are described in proof of Theorem~\ref{equi}.

\begin{corollary}
    Let $G$ be a strong digraph.
    Then, the univariate distance ideal $U_2(G) \neq \langle 1 \rangle$ if and only if $G$ is either $\overrightarrow{C}_3$ or $\Lambda(a,b,c,d)$, when 
    \begin{multicols}{3}
    \begin{itemize}
        \item $a\neq0,b=c=d=0$, 
        \item $c\neq0,a=b=d=0$, or
        \item $b,d\geq1, a=c=0$.
    \end{itemize}
    \end{multicols}
\end{corollary}
\begin{proof}
    It follows by evaluating $I_2(\overrightarrow{C}_3)$ and $I_2(\Lambda(a,b,c,d))$ at $\{x_i=t\}$ and $\{y_j=t\}$.
\end{proof}

If we evaluate $I_2(\overrightarrow{C}_3)$ and $I_2(\Lambda(a,b,c,d))$ at $\{x_i=0\}$ and $\{y_j=0\}$, up to graphs with only one vertex, the obtained ideal is generated by 1.
Therefore, the second invariant factor of the SNF of the distance matrix of any strong digraph with at least 2 vertices is equal to 1.

Now, we are going to give a description of the distance ideals of the strong digraphs described by $\Lambda(a,b,0,1)$.
Previously in \cite{at,corrval}, the distance ideals of complete graphs and star graphs were computed.
These graphs regarded as digraphs, are contained in the strong digraphs described by $\Lambda(a,b,0,1)$.
Thus, this extends previous results.
Let us recall the distance ideals of complete graphs and star graphs.

\begin{proposition}\cite{corrval}\label{teo:criticalidealscompletegraph}
    For $k\in\{1,\dots,n-1\}$, the set
    $
        \left\{ \prod_{i \in \mathcal{I}}(x_i+1) \, \big| \, \mathcal{I}\subset [n], \, |\mathcal{I}|=k-1 \right\}
    $
    is a reduced Gr\"obner basis of $I_k(K_n)$. 
    And
    $
        \det(D_X(K_n))=\prod_{i=1}^n(x_i+1)-\sum_{j=1}^n\prod_{i\neq j}(x_i+1)
    $.
\end{proposition}

\begin{proposition}\cite{at}
    For $k \in \{1,\dots,m\}$,  the distance ideal $I_k(K_{m,1})$
    is generated by 
    \[
    \left\{ (2y-1)\prod\limits_{i\in \mathcal{I}}(x_i-2) : \mathcal{I}\subset [m] \text{ and } |\mathcal{I}| = k-2  \right\}
    \bigcup
    \left\{ \prod\limits_{i\in \mathcal{I}}(x_i-2) : \mathcal{I}\subset [m] \text{ and } |\mathcal{I}| = k-1\right\}.
    \]
    And, the determinant of $D_X(K_{m,1})$ is equal to
    \begin{equation}
        y\prod\limits_{i=1}^m(x_i-2) + (2y-1)\sum\limits_{i=1}^m\prod\limits^m_{\substack{j=1\\ j\neq i}}(x_j-2).
    \end{equation}
\end{proposition}

In \cite[Theorem 3]{Nbook}, it was proven that if $\sf M$ and $\sf N$ are two equivalent matrices with entries in, then, for each $i\geq0$, the determinantal ideals generated by the $i$-minors of $\sf M$ and $\sf N$, respectively, coincide.
Therefore, the next result follows.

\begin{lemma}\label{lemma:equivLmatrix}
    Suppose $D_X(\Lambda(a,b,0,1))$ is 
    \[
    \diag(z_1,\dots,z_a,y_{1},\dots,y_{b},x)+
    \begin{bmatrix}
        {\sf J}_a - {\sf I}_a & {\sf J}_{a,b} & 2{\sf J}_{a,1}\\
        2{\sf J}_{b,a} & 2{\sf J}_b - 2{\sf I}_b & {\sf J}_{b,1}\\
        {\sf J}_{1,a} & {\sf J}_{1,b} & 0
    \end{bmatrix},
    \]
    where $\sf I$ and $\sf J$ denote identity matrix and the all-ones matrix, respectively.
    Therefore, $D_X(\Lambda(a,b,0,1))$ is equivalent to
    \[
    {\sf L}=
    \diag(x,y_{1}\dots,y_{b},z_1,\dots,z_a)+
    \begin{bmatrix}
        0 & {\sf J}_{1,b} & {\sf J}_{1,a}\\
        (1-2x){\sf J}_{b,1} & -2{\sf I}_b & {\sf 0}_{b,a}\\
        (2-x){\sf J}_{a,1} & {\sf 0}_{a,b} & -{\sf I}_a
    \end{bmatrix},
    \]
    and $I_k(\Lambda(a,b,0,1))$ coincides with the $k$-{\it th} determinantal ideal of matrix $\sf L$.
\end{lemma}
\begin{proof}
    It follows since
    $$
\begin{bmatrix}
    x & 1   & \dots & 1 & 1 & \dots & 1 \\
    1 & y_1 & \dots & 2 & 2 & \dots & 2 \\
\vdots & \vdots & \ddots & \vdots & \vdots & \ddots & \vdots \\
    1 & 2 & \dots & y_b & 2 & \dots & 2 \\
    2 & 1 & \dots & 1 & z_1 & \dots & 1 \\
\vdots & \vdots & \ddots & \vdots & \vdots & \ddots & \vdots \\
    2 & 1 & \dots & 1 & 1 & \dots & z_a \\
\end{bmatrix} \sim
\begin{bmatrix}
    x & 1   & \dots & 1 & 1 & \dots & 1 \\
    1-2x & y_1-2 & \dots & 0 & 0 & \dots & 0 \\
\vdots &\vdots & \ddots &\vdots & \vdots & \ddots & \vdots \\
    1-2x & 0 & \dots & y_b-2 & 0 & \dots & 0 \\
    2-x & 0 &\dots & 0 & z_1-1 & \dots & 0 \\
\vdots & \vdots  & \ddots & \vdots & \vdots & \ddots & \vdots \\
    2-x & 0 & \dots & 0 & 0 & \dots & z_a-1 \\
\end{bmatrix},
$$
and the fact that if $\sf M\sim N$, then, for each $i\geq0$, the ideals generated by the $i$-minors of $\sf M$ and $\sf N$, respectively, coincide, which was proved in~\cite[Theorem 3]{Nbook}.
\end{proof}

\begin{theorem}\label{teo:kdistideallambda}
    For $k<a+b+1$, the $k$-{\it th} distance ideal of $\Lambda(a,b,0,1)$ is
    \begin{equation}\label{kthideal}
\left\langle \prod_{s_{i}\leq b}\frac{y_{s_{i}}-2}{y_{r}-2}\prod_{b<s_{i}\leq b+a}(z_{s_{i}-b}-1), \prod_{s_{i}\leq b}(y_{s_{i}}-2)\prod_{b<s_{i}\leq b+a}\frac{z_{s_{i}-b}-1}{z_{t-b}-1} \,\bigg|\,
\begin{array}{c}
    r\in\lbrace s_{i}\mid s_{i}\leq b\rbrace,\\
    t\in\lbrace s_{i}\mid b<s_{i}\leq b+a\rbrace,\\
    s_{1}<\cdots < s_{k-1}
\end{array}
\right\rangle.
\end{equation}

    And the $(a+b+1)$-{\it th} distance ideal is generated by the polynomial
    \begin{equation}\label{determinant}
    \left(x-\sum_{r=1}^{b}\frac{1-2x}{y_{r}-2}-\sum_{s=1}^{a}\frac{2-x}{z_{s}-1}\right)\prod_{r=1}^{b}(y_{r}-2)\prod_{s=1}^{a}(z_{s}-1).
    \end{equation}
\end{theorem}

\begin{proof}
    Let $\Lambda=\Lambda(a,b,0,1)$.
    By Lemma~\ref{lemma:equivLmatrix}, we will focus on computing the minors of 
    \[
    {\sf L}=
    \diag(x,y_{1}\dots,y_{b},z_1,\dots,z_a)+
    \begin{bmatrix}
        0 & {\sf J}_{1,b} & {\sf J}_{1,a}\\
        (1-2x){\sf J}_{b,1} & -2{\sf I}_b & {\sf 0}_{b,a}\\
        (2-x){\sf J}_{a,1} & {\sf 0}_{a,b} & -{\sf I}_a
    \end{bmatrix}.
    \]
    For $ k \geq 1 $, let $ \mathcal{I} $ and $ \mathcal{J} $ be two subsets of $k$ indices from $ \{1, 2, \ldots, a+b+1\} $, and $ {\sf L}[\mathcal{I,J}] $ denote the submatrix of the matrix {\sf L} formed by choosing the rows indexed by $ \mathcal{I} $ and the columns indexed by $ \mathcal{J} $. 
    The $ k $-minor $\text{det}({\sf L}[\mathcal{I,J}])$ of ${\sf L}$ depends on whether $ 1 \in \mathcal{I} \cup \mathcal{J} $ or not, and on the intersection $ \mathcal{I} \cap \mathcal{J} $.

    If $1 \notin \mathcal{I} \cup \mathcal{J}$, then $\det({\sf L}[\mathcal{I,J}])$ is either the product of the elements of the diagonal of $\sf L$ indexed by $\mathcal{I}$ and $\mathcal{J}$, or $0$, depending on whether $\mathcal{I} = \mathcal{J}$ or not.
    
    If $1 \in \mathcal{J}$ and $1 \notin \mathcal{I}$, then one of the two polynomials $1-2x$ or $2-x$ multiplies the product of the diagonal elements of $\sf L$ indexed by $\mathcal{I}$ and $\mathcal{J}$, provided $\mathcal{J} \setminus \{1\} \subseteq \mathcal{I}$. 
    Notice that if this last condition does not hold, ${\sf L}[\mathcal{I,J}]$ contains an all-zeros column. 
    Finally, when $1\in \mathcal{I}\cap \mathcal{J}$, there are three possible scenarios based on whether the symmetric difference $\mathcal{I}\Delta \mathcal{J}$ contains more than two, exactly two or zero elements. 
    The first case leads with a zero determinant since an all-zero row appears. 
    The second case results in the product of the diagonal elements of $\sf L$ indexed by $\mathcal{I}$ and $\mathcal{J}$ multiplied by the product of two (not necessarily different) polynomials from $1-2x$ or $2-x$. 
    When $\mathcal{I=J}$, just the identity and permutations of the form $(1,j)$ will not vanish the product on the sum of the determinant. 
    Therefore, $\det({\sf L}[\mathcal{I,J}])$ equals to 
    $$
    \left(x-\sum_{r_{i}\leq b}\frac{1-2x}{y_{r_{i}}-2}-\sum_{b<r_{i}\leq b+a}\frac{2-x}{z_{r_{i}-b}-1}\right)\prod_{r_{i}\leq b}(y_{r_{i}}-2)\prod_{b<r_{i}\leq b+a}(z_{r_{i}-b}-1).$$
    Where $\mathcal{I=J}=\{1,r_{1},r_{2},\dots, r_{k-1}\}$.
    
    Therefore, for $k<a+b+1$, the $k-$th distance ideal $I_{k}(\Lambda(a,b,0,1))$ is 
    
    \begin{equation}\label{eq:ideal Lambda(a,b,0,1)}
    \left\langle \prod_{s_{i}\leq b}\frac{y_{s_{i}}-2}{y_{r}-2}\prod_{b<s_{i}\leq b+a}(z_{s_{i}-b}-1), \prod_{s_{i}\leq b}(y_{s_{i}}-2)\prod_{b<s_{i}\leq b+a}\frac{z_{s_{i}-b}-1}{z_{t-b}-1} \,\bigg|\,
    \begin{array}{c}
        r\in\lbrace s_{i}\mid s_{i}\leq b\rbrace,\\
        t\in\lbrace s_{i}\mid b<s_{i}\leq b+a\rbrace,\\
        s_{1}<\cdots < s_{k-1}
    \end{array}
    \right\rangle.
    \end{equation}
    
    And the $(a+b+1)$-{\it th} ideal is generated by the determinant of the matrix which is equal to
    
    \begin{equation}\label{eq:determinant Lambda(a,b,0,1)}
    \left(x-\sum_{r=1}^{b}\frac{1-2x}{y_{r}-2}-\sum_{s=1}^{a}\frac{2-x}{z_{s}-1}\right)\prod_{r=1}^{b}(y_{r}-2)\prod_{s=1}^{a}(z_{s}-1).
    \end{equation}
\end{proof}

The following is a consequence of previous theorem.

\begin{theorem}\label{teo:snfLambda(a,b,0,1)}
    The Smith Normal form of $D(\Lambda(a,b,0,1))$ is ${\sf I}_{a+2}\oplus 2{\sf I}_{b-2} \oplus [8a+2b]$ for $b\geq 2$ and $\textrm{SNF}(D(\Lambda (a,1,0,1)))={\sf I}_{a+1}\oplus [4a+1]$.
\end{theorem}

\begin{proof}
We can recover the invariant factors as described in Proposition~\ref{prop:eval1} by evaluating the variables in expressions~(\ref{kthideal}) and (\ref{determinant}) of Theorem~\ref{teo:kdistideallambda} at $0$. 

Notice that there are only $k-1$ factors in the products of generators depicted at the distance ideal described in~(\ref{kthideal}). 
Among this $k-1$ factors, at least $k-2$ can appear as variables $y_{s_{i}}-2$ due to the quotient $(y_{s_{i}}-2)/(y_{r}-2)$. 
By evaluating the variables at $0$, this proves that the first $a+2$ distance ideals are trivial. 
For $a+2<k\leq a+b$, the minimum number of $-2$ factors that appear in the ideal is $k-a-2$. 
Consequently, the $k$-{\it th} distance ideal is generated by $2^{k-a-2}$. 
The last ideal is generated by the expression in~(\ref{determinant}) evaluated at zero. 
A calculation shows that the determinant of the $(a+b+1)$-{\it th} distance matrix is $(-1)^{a+b}\left(b2^{b-1}+a2^{b+1}\right)$. 
According to Proposition~\ref{prop:eval1}, the invariant factors of the distance matrix of $\Lambda(a,b,0,1)$ are $1,1,\dots, 1,1,2,2,\dots, 2, 2b+8a$.

If $b=1$, then just exists one variable $y_{b}-2$ may appear in the product of polynomials of~\ref{eq:ideal Lambda(a,b,0,1)}. 
However, this factor disappears because of the quotient $(y_{b}-2)/(y_{b}-2)$. 
Therefore, $\Delta_{i}(D(\Lambda(a,1,0,1)))=1$ for $k=1,\dots,a+1$. 
Formula~(\ref{eq:determinant Lambda(a,b,0,1)}) holds even for $b=1$ evaluated at $0$ for the determinant of $D(\Lambda(a,1,0,1))$. 
This calculation yields with $\text{det}(D(\Lambda(a,1,0,1)))=(-1)^{a+1}(1+4a)$. 
Thus, the invariant factors of $D(\Lambda(a,1,0,1))$ are $1,1,\dots, 1+4a$. 
Finally, $\textrm{SNF}(D(\Lambda (a,1,0,1)))={\sf I}_{a+1}\oplus [4a+1]$.
\end{proof}

Now we give analogous results for $\Lambda(a,1,0,d)$.

\begin{lemma}\label{lemma:equivLmatrix2}
    Suppose $D_X(\Lambda(a,1,0,d))$ is 
    \[
    \diag(z_1,\dots,z_a,y,x_1\dots x_d)+
    \begin{bmatrix}
        {\sf J}_a - {\sf I}_a & {\sf J}_{a,1} & 2{\sf J}_{a,d}\\
        2{\sf J}_{1,a} & 0 & {\sf J}_{1,d}\\
        {\sf J}_{d,a} & {\sf J}_{d,1} & 2{\sf J}_d - 2{\sf I}_d
    \end{bmatrix},
    \]
    where $\sf I$ and $\sf J$ denote identity matrix and the all-ones matrix, respectively.
    Therefore, $D_X(\Lambda(a,1,0,d))$ is equivalent to
    \[
    {\sf L}=
    \diag(y,z_{1}\dots,y_{a},x_1,\dots,x_d)+
    \begin{bmatrix}
        0 & (2-y){\sf J}_{1,a} & (1-2y){\sf J}_{1,d}\\
        {\sf J}_{a,1} & -{\sf I}_a & {\sf 0}_{a,d}\\
        {\sf J}_{d,1} & {\sf 0}_{d,a} & -2{\sf I}_a
    \end{bmatrix},
    \]
    and $I_k(\Lambda(a,1,0,d))$ coincides with the $k$-{\it th} determinantal ideal of matrix $\sf L$.
\end{lemma}
\begin{proof}
    It follows since 
\[
\begin{bmatrix}
    y & 2 & \dots & 2 & 1 & \dots & 1 \\
    1 & z_1 & \dots & 1 & 2 & \dots & 2 \\
\vdots & \vdots & \ddots & \vdots & \vdots & \ddots & \vdots \\
    1 & 1 & \dots &z_{a} & 2 & \dots & 2 \\
    1 & 1 & \dots & 1 & x_1 & \dots & 2 \\
\vdots & \vdots & \ddots & \vdots & \vdots & \ddots & \vdots \\
    1 & 1 & \dots & 1 & 2 & \dots & x_{d} \\
\end{bmatrix} \sim
\begin{bmatrix}
    y & 2-y & \dots & 2-y & 1-2y & \dots & 1-2y \\
    1 & z_{1}-1 & \dots & 0 & 0 & \dots & 0 \\
\vdots &\vdots & \ddots &\vdots & \vdots & \ddots & \vdots \\
    1 & 0 & \dots & z_{a}-1 & 0 & \dots & 0 \\
    1 & 0 &\dots & 0 & x_{1}-2 & \dots & 0 \\
\vdots & \vdots  & \ddots & \vdots & \vdots & \ddots & \vdots \\
    1 & 0 & \dots & 0 & 0 & \dots & x_{d}-2 \\
\end{bmatrix}.
\]
\end{proof}

An analogous argument as in proof of Theorem~\ref{teo:kdistideallambda} will give us next result.

\begin{theorem}\label{teo:kdistideallambda2}
    For $k<a+d+1$, the $k$-{\it th} distance ideal of $\Lambda(a,1,0,d)$ is
    \begin{equation}
\left\langle \prod_{s_{i}\leq d}\frac{x_{s_{i}}-2}{x_{r}-2}\prod_{d<s_{i}\leq d+a}(z_{s_{i}-d}-1), \prod_{s_{i}\leq d}(x_{s_{i}}-2)\prod_{d<s_{i}\leq d+a}\frac{z_{s_{i}-d}-1}{z_{t-d}-1} \,\bigg|\,
\begin{array}{c}
    r\in\lbrace s_{i}\mid s_{i}\leq d\rbrace,\\
    t\in\lbrace s_{i}\mid d<s_{i}\leq d+a\rbrace,\\
    s_{1}<\cdots < s_{k-1}
\end{array}
\right\rangle.
\end{equation}
    And the $(a+b+1)$-{\it th} distance ideal is generated by the polynomial
    \begin{equation}\label{determinantLambda(a,1,0,d)}
\left(y-\sum_{r=1}^{d}\frac{1-2y}{x_{r}-2}-\sum_{s=1}^{a}\frac{2-y}{z_{s}-1}\right)\prod_{r=1}^{d}(x_{r}-2)\prod_{s=1}^{a}(z_{s}-1).
\end{equation}
\end{theorem}

\begin{corollary}
    The Smith Normal form of $D(\Lambda(a,1,0,d))$ is ${\sf I}_{a+2}\oplus 2{\sf I}_{d-2} \oplus [8a+2d]$ for $d\geq 2$.
\end{corollary}

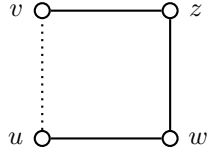
\begin{figure}[h]
    \centering
    \begin{tikzpicture}[scale=1.7]
        \node[label={180:$u$}] (z) at (0,0){};
        \node[label={0:$w$}] (w) at (1,0){};
        \node[label={0:$z$}] (y) at (1,1){};
        \node[label={180:$v$}] (x) at (0,1){};
        \draw (z) to (w);
        \draw (w) to (y);
        \draw (y) to (x);
        \draw[dotted] (z) to (x);
    \end{tikzpicture}
    \caption{Pattern $\mathcal{F}_6$}
    \label{fig:F0}
\end{figure}

Using the concept of pattern in the simple connected graphs setting, we obtain the following characterization of simple connected graphs with one trivial distance ideal.

\begin{proposition}
    Let $G$ be a connected simple graph.
    Then $G$ has exactly one trivial distance ideal if and only if pattern $\mathcal{F}_6$, described in Figure~\ref{fig:F0}, is not contained in $G$.    
\end{proposition}
\begin{proof}
    It follows by observing that $G$ does not contain pattern $\mathcal{F}_6$ if and only if $G$ is $\{P_4,\sf{paw},\sf{diamond}\}$-free.
\end{proof}

\section{Circuits}\label{sec:circuits}

In this section we explore distance ideals of circuits, which are not simple.
We start by computing the SNF of the distance matrix of circuits. 
From this, we observe that the third invariant factor of $\snf\left(D\left(\overrightarrow{C}_n\right)\right)$ is different of 1, which implies that $I_3\left(\overrightarrow{C_n}\right)$ is non-trivial, for $n\geq3$.
Then, we compute a basis for $I_3\left(\overrightarrow{C_n}\right)$.
The other distance ideals are more difficult, therefore, we only give some insights.




\begin{proposition}\label{prop:SNFcircuits}
The Smith normal form of the distance matrix of the circuit digraph $\overrightarrow{C_n}$ with $n\geq3$ is ${\sf I}_2\oplus n{\sf I}_{n-3}\oplus [n^2(n-1)/2)]$.
\end{proposition}
\begin{proof}
First note that the distance matrix $D(\overrightarrow{C_n})$ of the circuit $\overrightarrow{C_n}$ with $n$ vertices is
\[
\begin{bmatrix}
0 & 1 & \cdots & n-1 \\
n-1 & 0 & \cdots & n-2 \\
\vdots & \vdots & \ddots & \vdots \\
1 & 2 & \cdots & 0
\end{bmatrix}.
\]
To obtain the SNF, we will perform elementary operations on $D(\overrightarrow{C_n})$.
First, by subtracting the $(n-1)$-{\it th} column to the last column, and subtracting the $(n-2)$-{\it th} column to the $(n-1)$-{\it th} column, and so on, just leaving the first column the same, we obtain the matrix
\[
\begin{bmatrix}
0 & 1 & 1 & \cdots &1\\
n-1 & -n+1 & 1 & \cdots &1 \\
n-2 & 1 & -n+1& \cdots & 1\\
\vdots & \vdots & \vdots & \ddots & \vdots \\
1 & 1 & 1 & \cdots & -n+1
\end{bmatrix}.
\]
From this matrix, by subtracting the third column to the second, the fourth to the third, and so on, we obtain the matrix
\[
\begin{bmatrix}
0 & 0 & 0 & \cdots & 0 & 1\\
n-1 & -n & 0 & \cdots & 0 & 1 \\
n-2 & n & -n & \cdots & 0 & 1\\
\vdots & \vdots & \vdots & \ddots & \vdots& \vdots \\
2 & 0 & 0  & \cdots & -n & 1\\
1 & 0 & 0  & \cdots & n & -n+1
\end{bmatrix}
\]
Now, by adding the first column multiplied by $-n$ to the $(n-1)$-{\it th} column, we obtain the matrix
\[
\begin{bmatrix}
0 & 0 & 0 & \cdots & 0 & 1\\
n-1 & -n & 0 & \cdots & -(n-1)n & 1 \\
n-2 & n & -n & \cdots & -(n-2)n & 1\\
\vdots & \vdots & \vdots & \ddots & \vdots& \vdots \\
3 & 0 & 0  & \cdots & -3n & 1\\
2 & 0 & 0  & \cdots & -3n & 1\\
1 & 0 & 0  & \cdots & 0 & -n+1
\end{bmatrix}
\]
From this matrix, we can observe that with the first row, the elements, except for the first one, of the last column can be eliminated. 
And, with the last row, the elements, except for the last element, of the first column can be eliminated, yielding the matrix
\[
{\sf I}_2
\oplus
(-1)\begin{bmatrix}
n &0 & \cdots &0 & (n-1)n \\
-n &n & \cdots &0 & (n-2)n \\
\vdots & \vdots & \ddots & \vdots & \vdots \\
0 & 0 & \cdots & n & 3n \\
0 & 0 & \cdots & -n & 3n
\end{bmatrix}
\]
Then, the scalar multiplying the big matrix can be eliminated, and adding the first row of the big matrix to the second row, the third row to the fourth and so on, yields the matrix
\[
{\sf I}_2
\oplus
\begin{bmatrix}
n &0 & \cdots &0 & (n-1)n \\
0 &n & \cdots &0 & (n-1+n-2)n \\
\vdots & \vdots & \ddots & \vdots & \vdots \\
0 & 0 & \cdots & n & (n-1+n-2+\cdots+3)n \\
0 & 0 & \cdots & 0 & (n-1+n-2+\cdots+3+3)n
\end{bmatrix}
\]
All the elements of the last column are multiple of $n$, then they can be eliminated with the previous columns, except for the last one, yielding the matrix ${\sf I}_2\oplus n{\sf I}_{n-3}\oplus [n^2(n-1)/2)]$.
\end{proof}

A consequence of Propositions~\ref{prop:eval1} and \ref{prop:SNFcircuits} is that the third distance ideal of circuits with at least 3 vertices is not trivial.
Therefore, it is natural to ask for a basis for such ideal.
This is what we will compute next.
First note that the second distance ideal of $\overrightarrow{C}_n$ with $n\geq 4$ is trivial since the matrix $D_X(\overrightarrow{C}_n)$ contains the submatrix
\[
\begin{bmatrix}
    2 & 3\\
    1 & 2 
\end{bmatrix}.
\]
The third distance ideal of $\overrightarrow{C}_n$ with $n\geq 4$ is more interesting.
First, for $n=4$, $I_3\left(\overrightarrow{C}_4\right)$ is generated by the $\binom{4}{3}^2$ polynomials, which is non-trivial since the third invariant factor of the SNF of $D(\overrightarrow{C}_4)$ is different of 1.
The reader might compute a Gr\"oebner basis of this ideal with the \texttt{sagemath} code of Appendix~\ref{appendix:distanceideals}.
For $n\geq 5$, we are able to give a simple description of the third distance ideal, which we give next.

\begin{theorem}\label{coro:teo:3idealciclo}
    The third distance ideal of $\overrightarrow{C}_n$, with $n\geq 5$, is $\langle x_1,x_2,\ldots ,x_n,n\rangle$.
\end{theorem}

\begin{proof}
    Let ${\sf M} = D_X(\overrightarrow{C}_n)$.
    Let $\mathcal{I,J}\subset [n]$ and let ${\sf M}[\mathcal{I, J}]$ be the submatrix of $\sf M$ obtained by choosing the indices in $\mathcal{I}$ and $\mathcal{J}$ as the rows and columns of the submatrix respectively. 
    Now, let $\mathcal{I}=\{i_1,i_2,i_3\}$ and $\mathcal{J}=\{j_1,j_2,j_3 \}$ with $1\leq i_1 <i_2<i_3\leq n$ and with $1\leq j_1 <j_2<j_3\leq n$. Then we have three types of such submatrices. 
    The first type occurs when $|\mathcal{I\cap J}|=\emptyset$. 
    In this case, we have
    \[
    \det({\sf M}[\mathcal{I,J}])=
    \det
    \begin{bmatrix}
    n-i_1+j_1 & n-i_1+j_2 & n-i_1+j_3 \\
    n-i_2+j_1 & n-i_2+j_2 & n-i_2+j_3\\
    n-i_3+j_1 & n-i_3+j_2 & n-i_3+j_3
    \end{bmatrix}
    \]

    Note that the numbers on the entries of the previous matrix should be taken modulo $n$. 
    If we compute the determinant as presented we have that $\det({\sf M}[\mathcal{I,J}])=0$, therefore taking the determinant expression with each entry $\pmod n$, we have that $\det({\sf M}[\mathcal{I,J}])\equiv 0 \mod n$. Thus $\det({\sf M}[\mathcal{I,J}])=0$ or a multiple of $n$.

    Now, for the case $|\mathcal{I\cap J}|=3$. 
    Let us set $\mathcal{I=J}=\{i,j,k\}$, and without loss of generality, we assume that $i<j<k$. 
    Then we have

    \[
        \det({\sf M}[\mathcal{I,J}])=\det\begin{bmatrix}
        x_i & -i+j & -i+k \\
        n-j+i & x_j & -j+k\\
        n-k+i & n-k+j & x_k
        \end{bmatrix} 
    \]
    That is 
    \begin{align*}
        \det({\sf M}[\mathcal{I,J}])&=x_ix_jx_k-x_i((n-k+j)(-j+k))-x_j((n-k+i)(-i+k))\\
        &-x_k((n-j+i)(-i+j) ) + n(-j^2+j(i+k)-k^2+(k-i)(i+n)).
    \end{align*}
    Now let $\mathcal{I\cap J}=\{i,j\}$, $i<j$ and let also $\mathcal{I}=\{i,j,k_1\}$ and $\mathcal{J}=\{i,j,k_2\}$. 
    Moreover, let $k_1<k_2$. 

In this case, we have
\[
\det({\sf M}[\mathcal{I,J}])=
\det\begin{bmatrix}
    x_i & -i+j & (n-i+k_2)\mod n \\
    n-j+i & x_j & (n-j+k_2)\mod n\\
    (n-k_1+i)\mod n & (n-k_1+j)\mod n & -k_1+k_2
\end{bmatrix}
\]

That is,
\begin{align*}
    \det({\sf M}[\mathcal{I,J}])&=x_ix_j(-k_1+k_2)-x_i((n-k_1+j)_n(n-j+k_2)_n)-x_j((n-k_1+i)_n(n-i+k_2)_n)\\
    &-(-k_1+k_2)((n-j+i)(-i+j) ) + n(-j^2+j(i+k)-k^2+(k-i)(i+n)) +cn 
\end{align*}

Finally, let $\mathcal{I\cap J}=\{ k \}$, and let also $\mathcal{I}=\{i_1,i_2,k\}$ with $i_1<i_2$ and $\mathcal{J}=\{j_1,j_2,k\}$ with $j_1<j_2$. Without loss of generality, also set $i_1<j_1$. 

In this case, we have
\[
\det({\sf M}[\mathcal{I,J}])=
\det\begin{bmatrix}
    x_k & (n-k+j_1)\mod n & (n-k+j_2) \mod n \\
    (n-i_1+k)\mod n & (n-i_1+j_1)\mod n & (n-i_1+j_2)\mod n\\
    (n-i_2+k)\mod n & (n-i_2+j_1)\mod n & (n-i_2+j_2)\mod n
\end{bmatrix}
\]

That is,
\begin{align*}
    \det({\sf M}[\mathcal{I,J}]) = & x_k((n-i_1+j_1)\mod n(n-i_2+j_2)\mod n)\\
    & -x_k((n-i_2+j_1)\mod n(n-i_1+j_2)\mod n)\\
    & +cn(i_1-i_2)(j_2-j_1)
\end{align*}

In particular, given $1\leq k\leq n$ and $n\geq 5$, we can choose $\mathcal{I}=\{(k-2)\mod n,(k-1)\mod n,k\mod n\}$ and $\mathcal{J}=\{k\mod n,(k+1)\mod n,(k+1)\mod n\}$ and the corresponding polynomial in the 3-{\it rd} distance ideal is simply $x_k$. 
Moreover, for any such $k$, we have that $x_k+n\in I_3(\overrightarrow{C_n})$, where $n\geq 5$. 
We can show this claim by choosing $\mathcal{J}=\{(k-2)\mod n,(k-1)\mod n,k\mod n\}$ and $\mathcal{I}=\{(k)\mod n,(k+1)\mod n,(k+1)\mod n\}$. 
Then, also $n\in I_3(\overrightarrow{C_n})$.

Thus, a base for the third ideal is given by $\{x_1,\ldots,x_n,n\}$ since the constant monomial of any polynomial in the ideal is a multiple of $n$.
\end{proof}

A consequence of Theorem~\ref{coro:teo:3idealciclo} is the following.

\begin{corollary}
    For $n\geq5$, $\overrightarrow{C_n}$ is in the family of strong digraphs with exactly 2 trivial distance ideals.
\end{corollary}

This also implies that there is an infinite number of minimal forbidden induced subdigraphs for the family of strong digraphs with exactly one trivial distance ideals.
Therefore, obtaining a classification theorem for the strong digraphs in $\Gamma_{1}$ in terms of induced subdigraphs seems more difficult than the obtained by using patterns as in Theorem~\ref{equi}.

It will be great to find a basis for the rest of the distance ideals of the circuits, but seems to be a hard problem, even for univariate distance ideals.
Since the univariate distance ideals can be obtained by evaluating the variables at $x_i=t$, then the first three univariate distance ideals of circuits are known.
For the other univariate distance ideal, we have the following conjecture which is based on computational experiments.

\begin{conjecture}\label{conj:univdistidealsC_n}
    Let $\overrightarrow{C_n}$ be the circuit with $n$ vertices and let $\omega = e^{2i\pi/n}$.
    The univariate distance ideals $U_k(\overrightarrow{C_n})$ are $\left\langle \bigcup_{i=0}^{k-2} t^{k-2-i}n^{i}\right\rangle$, for $4\leq k\leq n-2$.
    Moreover, $p(t)=\det(D_t(\overrightarrow{C_n}))=\prod_{j=0}^{n-1} (t+\omega^j+2\omega^{2j}+\cdots +(n-1)\omega^{(n-1)j})$.
\end{conjecture}

In general, it seems harder to find an expression for the base of the $(n-1)$-{\it th} univariate distance ideal of the circuits. 

\section*{Acknowledgement}
The authors were supported by SNI and CONACyT.
The research of Ralihe R. Villagran was supported, in part, by the National Science Foundation through the DMS Award $\#$1808376 which is gratefully acknowledged.

\appendix

\section{Computing distance ideals with \texttt{sagemath}}\label{appendix:distanceideals}
The following {\tt sagemath} code computes the distance ideals of the complete digraph with 3 vertices.
\begin{lstlisting}
def DistanceIdealsZZ(G):
    n = G.order()
    S = "["
    for i in range(n):
        if i > 0 :
            S += ","
        S += "x" + str(i)
    S += "]"
    R = PolynomialRing(ZZ, "x", n);
    DX = diagonal_matrix(list(R.gens())) + G.distance_matrix()
    print("D_X matrix:")
    print(DX)
    for i in range(1, n + 1):
        I = R.ideal(DX.minors(i)).groebner_basis()
        print("Distance ideals of size " + str(i) + ":")
        print(I)
g = digraphs.Complete(3)
g.show()
DistanceIdealsZZ(g)
\end{lstlisting}
The following {\tt sagemath} code computes the univariate distance ideals of the circuit digraph with 3 vertices.
Note we use two variables to be able to compute Gröbner bases, since the algorithm is not implemented over univariate polynomial rings.
\begin{lstlisting}
def UnivariateDistanceIdealZZ(G):
    n = G.order()
    S = ["t","s"]
    R = PolynomialRing(ZZ, S);
    R.inject_variables(verbose=False)
    Dt = t*identity_matrix(n) + G.distance_matrix()
    print("D_t matrix:")
    print(Dt)
    for i in range(1, n + 1):
        I = R.ideal(Dt.minors(i)).groebner_basis()
        print("Univariate distance ideals of size " + str(i) + ":")
        print(I)
G = digraphs.Circuit(3)
G.show()
UnivariateDistanceIdealZZ(G)
\end{lstlisting}
In particular, the following code computes the second distance ideal of the digraphs $\Lambda$, defined in Theorem~\ref{equi}.
\begin{lstlisting}
def SecondDistanceIdealZZ(G):
    n = G.order()
    S = "["
    for i in range(n):
        if i > 0 :
            S += ","
        S += "x" + str(i)
    S += "]"
    R = PolynomialRing(ZZ, "x", n);
    DX = diagonal_matrix(list(R.gens())) - G.distance_matrix()
    Minors = DX.minors(2)
    I = R.ideal(Minors).groebner_basis()
    print("\nGrobner Basis:")
    print(I)

def Lambda(a,b,c,d):
    print("a="+str(a)+"\nb="+str(b)+"\nc="+str(c)+"\nd="+str(d))
    G = digraphs.Complete(a)
    G.add_vertices([a+i for i in range(b)])
    G.add_edges([(i,a+j) for i in range(a) for j in range(b)])
    G.add_vertices([a+b+i for i in range(c)])
    G.add_edges([(a+i,a+b+j) for i in range(b) for j in range(c)])
    G.add_edges([(a+b+i, a+b+j) for i in range(c) for j in range(c) if i != j])
    
    G.add_vertices([a+b+c+i for i in range(d)])
    G.add_edges([(a+b+i,a+b+c+j) for i in range(c) for j in range(d)])
    G.add_edges([(a+i,a+b+c+j) for i in range(b) for j in range(d)])
    G.add_edges([(a+b+c+i, a+j) for i in range(d) for j in range(b)])
    G.add_edges([(a+b+c+i, j) for i in range(d) for j in range(a)])
    
    return G

a = 1; b = 1; c = 3; d = 1
SecondDistanceIdealZZ(Lambda(a,b,c,d))
\end{lstlisting}

\section{Varieties}\label{appendix:varieties}

The following {\tt sagemath} code draws a partial view of the variety associated with the third distance ideal of the complete digraph with 3 vertices.
\begin{lstlisting}
P.<x,y,z> = PolynomialRing(ZZ, 3, order = "lex")
var("x,y,z")
g = DiGraph({0:[1,2],1:[0,2],2:[0,1]})
g.plot()
d = g.distance_matrix()+diagonal_matrix([x,y,z])
print(d.det())
l = 400
u = 400
implicit_plot3d(d.det() == 0, (x,-l,u), (y,-l,u), (z,-l,u), color = "lightcoral", plot_points = 100, frame = False)   
\end{lstlisting}
The following {\tt sagemath} code draws a partial view of the variety associated with the third distance ideal of the circuit digraph with 3 vertices.
\begin{lstlisting}
P.<x,y,z> = PolynomialRing(ZZ, 3, order = "lex")
var("x,y,z")
g = DiGraph({0:[1],1:[2],2:[0]})
g.plot()
d = g.distance_matrix() + diagonal_matrix([x,y,z])
print(d.det())
l = 3
u = 10
implicit_plot3d(d.det() == 0, (x,-l,u), (y,-l,u), (z,-l,u), color = "lightcoral", plot_points = 100, frame = False)
\end{lstlisting}

\section{Strong digraphs with 4 vertices and non-trivial second ideal}\label{appendix:digraphssecondnontrivialideal}

\begin{lstlisting}
def SecondDistanceIdealZZ(G):
    n = G.order()
    S = "["
    for i in range(n):
        if i > 0 :
            S += ","
        S += "x" + str(i)
    S += "]"
    R = PolynomialRing(ZZ, "x", n);
    DX = diagonal_matrix(list(R.gens())) + G.distance_matrix()
    I = R.ideal(DX.minors(2)).groebner_basis()
    return I

for n in range(2,5):
    for g in digraphs(n):
        if g.is_strongly_connected():
            si = SecondDistanceIdealZZ(g)
            if si[0] != 1 and si[0] != -1:
                g.show()
\end{lstlisting}

\end{document}